\documentclass[11pt,centertags,oneside]{amsart}
\usepackage{amsmath,amstext,amsthm,amscd,typearea}
\usepackage{amssymb}
\usepackage{a4wide}
\usepackage[mathscr]{eucal}
\usepackage{mathrsfs}
\usepackage{typearea}
\usepackage{charter}
\usepackage{pdfsync}
\usepackage{tikz}
\usepackage{physics}
\usepackage{mathtools}
\usepackage{fancyhdr}
\usepackage{verbatim}
\usepackage{subfiles}
\usepackage[normalem]{ulem}
\usetikzlibrary{decorations.pathreplacing}
\usepackage[a4paper,width=16.5cm,top=3cm,bottom=3cm]{geometry}

\numberwithin{equation}{section}

\allowdisplaybreaks
\emergencystretch=\maxdimen
\usepackage{multicol}
\usepackage{xcolor}
\usepackage[hidelinks]{hyperref}

\newtheorem{theorem}{Theorem}[section]
\newtheorem{definition}[theorem]{Definition}
\newtheorem{proposition}[theorem]{Proposition}
\newtheorem{corollary}[theorem]{Corollary}
\newtheorem{lemma}[theorem]{Lemma}
\newtheorem{remark}[theorem]{Remark}

\makeatletter
\newcommand*{\rom}[1]{\expandafter\@slowromancap\romannumeral #1@}
\makeatother

\newcommand{\mb}{\mathbb}

\newcommand{\mcr}{\mathscr}

\newcommand{\ovl}{\overline}

\newcommand{\bkt}[1]{\left({#1}\right)}
\newcommand{\set}[1]{\left\{{#1}\right\}}

\newtheoremstyle{prob}{3em}{3em}{}{0pt}{\bfseries}{.}{5pt plus 1pt minus 1pt}{\thmname{#1}\thmnumber{#2}}
\theoremstyle{prob}

\DeclareMathOperator{\supp}{supp}

\definecolor{han}{rgb}{1.0, 0, 0}

\newcommand{\rem}[1]{ }
\begin{document}
\title{Size of Exceptional Sets in Weakly Mixing Systems}

\author{Jiyun Park}
\address{Department of Mathematics, Stanford University}
\email{jiyunp@stanford.edu}

\author{Kangrae Park}
\address{Department of Mathematical Sciences, Seoul National University}
\email{kangrae.park@snu.ac.kr}

\date{\today}

\begin{abstract}
We study exceptional sets for the Chacon transformation and, more generally, for a class of cutting-and-stacking transformations called restrictive tight maps. For these systems we explicitly construct a universal exceptional set \(J\subseteq\mathbb N\), valid uniformly for all measurable pairs \(A,B\in\mcr B\), such that for every increasing function \(h:\mathbb N\to\mathbb R_{>0}\) diverging to infinity,
\[
\bigl|J\cap[0,n]\bigr|\le(\log n)^{h(n)}
\quad\text{for all sufficiently large }n.
\]
The Chacon transformation considered in this paper belongs to this class, giving a logarithmic-scale universal exceptional set for Chacon. We also prove that this logarithmic scale is essentially sharp: for every tight map with no spacers above the last subcolumn, i.e. \(s_{m-1}=0\), and every \(t>0\), there exist measurable sets \(A,B\) such that every exceptional set \(J\) for \((A,B)\), if one exists, satisfies
\[
\bigl|J\cap[0,n]\bigr|\ge(\log n)^t
\quad\text{for all sufficiently large }n.
\]
The construction is based on recursive formulas for return-time distributions arising from the cutting-and-stacking structure. As a complementary quantitative principle, we show that if the corresponding Ces\`aro weak-mixing averages satisfy a rate \(o(b_N)\), then \(J_{A,B}\) may be chosen so that
\[
\bigl|J_{A,B}\cap[0,N]\bigr|=o(Nb_N).
\]
We apply this quantitative implication from Ces\`aro rates to exceptional-set bounds to several weakly mixing models, including interval exchange transformations, translation flows, and substitution dynamical systems, under the regularity assumptions of the available quantitative estimates. We also construct a separate weakly mixing one-spacer rank-one example for which exceptional sets for a suitable pair have polynomial lower growth.
\end{abstract}

\maketitle

\section{Introduction}
\label{sec:intro}

A measure-preserving system \((X,\mcr{B},\mu,T)\) is said to be \emph{mixing} if
\begin{equation}
 \mu\bigl(A\cap T^{-n}B\bigr)\longrightarrow\mu(A)\,\mu(B)
 \quad(n\to\infty)
 \label{eq:strong-mixing}
\end{equation}
for all measurable sets \(A,B\subseteq X\). Mixing appears in many contexts—geodesic flows on negatively curved manifolds, Anosov flows on tori, and beyond—and underlies equidistribution, rapid decay of correlations, and statistical limit laws.

Many natural systems (for instance, typical non-rotation interval exchange transformations or the classical Chacon transformation) fail to be mixing yet still display an averaged form of decorrelation known as \emph{weak mixing}. In the discrete-time setting, this means that, for all measurable sets \(A,B\subseteq X\),
\[
C_N \,=\;
\frac{1}{N}\sum_{n=0}^{N-1}\bigl|\mu(A\cap T^{-n}B)-\mu(A)\,\mu(B)\bigr|
\;\longrightarrow\;0
\quad(N\to\infty).
\]
An alternate description of weak mixing for a transformation is the following, which dates back to Halmos \cite{Halmos}.

A transformation $T$ is weak mixing if and only if for every $A,B\in \mcr{B}$, there is a zero-density set $J_{A,B}\subseteq \mathbb{N}$ for which $$\mu(A\cap T^{-n}B) \to \mu(A)\mu(B)$$ as $n\to \infty$ with $n\notin J_{A,B}$ \cite[Theorem~2.36]{ward}.

Throughout this paper, we will refer to \(J_{A,B}\) as an \emph{exceptional set} of \(T\). That is, \(J_{A,B}\subseteq\mathbb N\) is exceptional for \((A,B)\) exactly when 
\[
\mu(A\cap T^{-n}B)\to\mu(A)\,\mu(B)
\quad(n\to\infty,\;n\notin J_{A,B}).
\]
No density or size restriction is included in the word ``exceptional'' unless it is explicitly stated.
Similarly, if \(\mcr C\subseteq\mcr B\) is any collection of measurable sets, we say \(J_{\mcr C}\) is exceptional for \(\mcr C\) when it is exceptional for every pair \((A,B)\in\mcr C\times\mcr C\). 

In some sense, the size of exceptional sets can be used to quantify how close a transformation or flow is to being mixing. For a fixed pair $(A,B)$, a smaller exceptional set indicates stronger decorrelation for that pair. For all pairs, the extreme case in which the exceptional set can be taken empty is strong mixing. Allowing density-zero exceptional sets gives weak mixing; in the countably generated case, this can be formulated using a universal density-zero exceptional set. The size of exceptional sets is therefore not merely a reformulation of weak mixing: as the polynomial lower-bound example in Subsection~\ref{subsec:rank-one-lower-bound} shows, weak mixing alone can coexist with polynomially large necessary exceptional sets. This motivates the following fundamental questions:
\begin{enumerate}
 \item How does a quantitative weak-mixing rate control the possible size of $J_{A,B}$?
 \item What dynamical information is reflected by the optimal size of exceptional sets?
 \item Given a measure-preserving action, is there a way to construct explicit exceptional sets for it?
\end{enumerate}

In this paper, the constructive direction is our main focus. We construct logarithmic-scale universal exceptional sets for restrictive tight maps, a cutting-and-stacking class that includes the Chacon transformation considered below, and thereby identify systems for which exceptional sets are much smaller than what one obtains from weak-mixing rates alone. We also prove a complementary quantitative principle showing that a weak-mixing rate gives a corresponding upper bound on the size of exceptional sets, and we apply this principle to several standard weakly mixing models. The precise statements are in the next section.

\subsection{Main results}
\label{subsec:main-results}

The following proposition shows that the size of exceptional sets can be bounded by the rate at which the relevant Ces\`aro averages converge.

\begin{proposition}
\label{prop:weak-conv-rate}
Fix $p\ge1$ and let $b_N \to 0$ be a positive sequence. Let $(X,\mcr{B},\mu,T)$ be a measure-preserving transformation and let $A,B\in\mcr{B}$ satisfy
 \[
 \frac{1}{N}\sum_{n=0}^{N-1}\bigl|\mu(A\cap T^{-n}B)-\mu(A)\mu(B)\bigr|^{p}
 =o(b_N)
 \quad(N\to\infty).
 \]
 Then there exists an exceptional set $J_{A,B}\subseteq\mathbb{N}$ such that
 \[
 \bigl|J_{A,B}\cap[0,N]\bigr|
 =o(N\,b_N)
 \quad\text{and}\quad
 \mu(A\cap T^{-n}B)\to\mu(A)\mu(B)
 \]
 as $n\to\infty$ with $n\notin J_{A,B}$. 
\end{proposition}

This result aligns with our intuition: if the Ces\`aro averages converge rapidly, only a few times can deviate significantly. Indeed, the proof is a straightforward quantitative refinement of the classical argument for \cite[Theorem~2.36]{ward} and is completed in Section~\ref{sec:exceptional-set}. There is no direct converse in this pairwise form: mixing systems may have $J_{A,B}=\varnothing$, while their quantitative correlation decay can still vary substantially. We apply this proposition in Section~\ref{sec:applications} to obtain exceptional-set bounds for several weakly mixing systems, under the regularity hypotheses of the cited quantitative estimates.

Our main focus is the explicit construction of exceptional sets, and hence the identification of systems for which exceptional sets are much smaller than what one obtains from weak-mixing rates alone. We work with a broad class of rank-one systems which we call \emph{restrictive tight maps} (Definitions~\ref{def1} and~\ref{def2}); the universal upper bound is proved for restrictive tight maps.

\begin{theorem}
\label{thm:general}
Let $(X, \mcr{B}, \mu, T)$ be a restrictive tight map. For any increasing $h : \mathbb{N} \to \mathbb{R}_{>0}$ diverging to infinity, there exists a set $J \subseteq \mathbb{N}$ such that $J$ is exceptional for $\mcr{B}$ and
\[
|J \cap [0, n]| \le (\log n)^{h(n)}
\]
for all sufficiently large $n$.
\end{theorem}

The Chacon transformation considered in this paper is a restrictive tight map. By Proposition~\ref{prop:chacon-prop}, it is weak mixing but non-mixing, hence it gives an iconic special case of the general result.

\begin{corollary}
\label{cor:chacon}
Let $T$ be the Chacon transformation. For any increasing $h : \mathbb{N} \to \mathbb{R}_{>0}$ diverging to infinity, there exists a set $J \subseteq \mathbb{N}$ such that
\[
|J \cap [0, n]| \le (\log n)^{h(n)}
\]
for all sufficiently large $n$, and $J$ is exceptional for every pair of Lebesgue-measurable sets $A,B\subseteq[0,1)$.
\end{corollary}

In fact, our construction yields exceptional sets for all $L^{2}$ functions (see Proposition~\ref{prop:J_f}). We remark that, after the first version of this paper was posted on arXiv, Moll \cite{moll2023speed} proved a quantitative weak-mixing estimate for the Chacon transformation for zero-mean Lipschitz observables tested against $L^2$ observables. Combined with Proposition~\ref{prop:weak-conv-rate}, this gives, for every $\gamma<1/6$, an automatic exceptional-set bound of order $o(n[\log_3 n]^{-\gamma})$ for the corresponding observable pairs. By choosing $h$ to grow sufficiently slowly, our bound for a universal exceptional set of the Chacon transformation is much stronger. Thus, for Chacon, direct construction gives much smaller exceptional sets than the automatic bounds obtained from the available weak-mixing-rate estimate alone. Furthermore, we show that the upper bounds above are essentially sharp:

\begin{theorem}
\label{thm:lower-bound}
Let $(X,\mcr B,\mu,T)$ be a tight map with spacer sequence $(s_0,\ldots,s_{m-1})$, and assume that no spacers are placed above the last subcolumn, i.e. $s_{m-1}=0$. For every $t>0$ there exist $A,B\in\mcr B$ such that if $J_{A,B}$ is any exceptional set for $A$ and $B$, then for some $N\in\mathbb{N}$,
\[
\lvert J_{A,B}\cap[0,n]\rvert\ge(\log n)^{t}\qquad(n\ge N).
\]
\end{theorem}

The logarithmic scale in Theorem~\ref{thm:general} should not be viewed as a consequence of weak mixing alone. In Subsection~\ref{subsec:rank-one-lower-bound} we construct a weakly mixing one-spacer rank-one transformation with polynomial lower bounds for exceptional sets: for every $\delta\in(0,1)$, there are a set $A$ and a constant $C>0$ such that every exceptional set $J$ for $(A,A)$ satisfies
\[
 |J\cap[0,N]|\ge C N^\delta
\]
for all sufficiently large $N$. Thus the logarithmic universal bound above reflects additional structure of restrictive tight maps, not weak mixing alone.

Note that, for a restrictive tight map, Theorem~\ref{thm:general} gives a \emph{universal} exceptional set; that is, $J$ is exceptional for all pairs of measurable sets in $\mcr B$. This naturally connects to the following proposition, which is essentially Corollary~3.2 of \cite{Friedman} (and so we do not repeat the proof).

\begin{proposition}
 \label{prop:weak-exceptional}
 Let $(X, \mcr{B}, \mu, T)$ be a measure-preserving system. If $\mcr{B}$ is countably generated, then $T$ is weak mixing if and only if there exists a zero-density set $J \subseteq \mathbb{N}$ that is exceptional for $\mcr{B}$.
\end{proposition}
\subsection{Related work}
The existence of weakly mixing systems that are not strongly mixing was shown by R. V. Chacon \cite{chacon} (see also \cite{chacon2}) and is commonly referred to as the Chacon transformation\footnote{There are more than one transformations that go by the same name of ``Chacon transformation'', and in particular the one we use here is not the one that appears in Chacon's original paper. Our definition of the Chacon map is given in Section~\ref{sec:chacon}.}. There are also other intermediate notions of mixing, such as mildly mixing (equivalently, having no nontrivial rigid factor; in particular $\liminf_{n\to\infty} \mu(A \triangle T^{-n} A) > 0$ for every $A$ with $0<\mu(A)<1$) and lightly mixing (i.e. $\liminf_{n\to\infty} \mu(A \cap T^{-n} B) > 0$ for all $A,B$ with $\mu(A),\mu(B)>0$), which are related in the following manner.
\[
\text{mixing}\;\Longrightarrow\;\text{lightly mixing}
\;\Longrightarrow\;\text{mildly mixing}
\;\Longrightarrow\;\text{weakly mixing}.
\]
The present paper takes a complementary quantitative viewpoint: rather than
placing a system in the mild/light mixing hierarchy, we measure how large the
set of times excluded from mixing-type convergence must be. We do not pursue a
general characterization of this size in terms of mild or light mixing.

The original Chacon map \cite{chacon} is known to be mildly mixing but not lightly mixing, while the map considered in our paper is known to be lightly mixing (but not strongly mixing) \cite{Friedman,Alpha}. The original Chacon map also has a trivial centralizer and minimal self-joinings of all orders. It is non-rigid and not isomorphic to its inverse \cite{Junco} \cite{Junco2} \cite{Fieldsteel} (see also \cite{Bell}). A transformation \(T\) is \(\alpha\)-mixing if there is a subsequence \(T^{m_k}\) converging weakly to \(\alpha\,\Theta + (1-\alpha)\mathrm{Id}\), where \(\Theta\) is the orthoprojector onto the constants. The Chacon transformation is not \(\alpha\)-mixing for \(0\le\alpha\le1\) \cite{Alpha}. King \cite{king1988joining} introduced the joining-rank invariant and analyzed the fine structural hierarchy of finite-rank mixing maps.

There are some generalizations of the Chacon transformation. T. Adams, N. Friedman, and C. Silva constructed an infinite measure-preserving rank-one transformation which can be viewed as a Chacon transformation in infinite measure \cite{Adams}. There is another version of an infinite Chacon transformation which has similar properties to the classical Chacon transformation \cite{Near}. V. V. Ryzhikov \cite{Ryz} generalized the Chacon transformation by using different sizes of spacers. The Chacon transformation is rank-one, and the examples and properties of rank-one transformations are studied in \cite{Ryz2}. A. del Junco and K. K. Park \cite{del1982example} constructed the first example of a measure-preserving flow with minimal self-joinings. There is a theory of joinings for two-dimensional Chacon-like transformations, establishing properties such as minimal self-joinings and disjointness \cite{park1991joinings}.

The Chacon transformation is an example of a rank-one transformation, a broader class of transformations that exhibit similar structural properties. Rank-one transformations have been extensively studied in various contexts, including mixing properties and factorization behaviors. Friedman and Ornstein \cite{friedman1972mixing} showed the existence of a mixing rank-one transformation. Moreover, rank-one mixing transformations are known to be mixing of all orders \cite{kalikow1984twofold, ryzhikov1993joinings}. Creutz and Silva \cite{creutz2010mixing} established that mixing in rank-one transformations is equivalent to the spacer sequence being slice-ergodic. The factorization properties of rank-one transformations have been extensively studied, particularly in connection with odometers and finite cyclic permutations \cite{foreman2023rank}. Spectral aspects of rank-one maps have also been investigated in relation to Mahler measure \cite{abdalaoui2021mahler}. Further results on rank-one transformations can be found in \cite{correia2024rank, creutz2024word, creutz2023measure}.

Rank-one transformations also serve as fundamental examples in the study of ergodic flows and higher-dimensional systems. There exists a rank-one infinite measure-preserving flow where every non-zero transformation possesses infinite ergodic index \cite{danilenko2011rank}. In the setting of the Chacon \(\mathbb{Z} \times \mathbb{Z}\) system, it has been shown that the time-zero partition is a generating partition under the transformation \(T\), and the centralizer of \(T\) is characterized \cite{johnson1997dynamical}.

Apart from the Chacon transformation, A. Katok \cite{katok1980interval} proved that interval exchange transformations (IETs) cannot be strong mixing, and A. Avila and G. Forni \cite{IET} proved that almost every typical (non-rotation) irreducible IET is weak mixing. If an IET can be obtained as an induced map of a rigid motion, then it is said to be of rotation class. For observables in the relevant regularity classes, Avila--Forni--Safaee \cite{IET_rate} obtain polynomial Ces\`aro correlation decay for typical non-rotation IETs and logarithmic decay in the rotation class. Quantitative weak-mixing and spectral estimates in related models, including translation flows and substitution systems, can be found in \cite{Forni2022,BufetovSolomyak2021,bufetov2025local,marshall2024quantitative,trevino2020quantitative}.

There are multiple properties concerning weak mixing rate. We say a weakly mixing transformation is \emph{partially weakly $f$-mixing} if there exists $A$ such that $C_N=o(f(N))$ for every $B$. A measure $\mu$ on a circle is called \emph{uniformly $f$-continuous} if there exists $C>0$ such that $\mu(I)\le C f(\abs{I})$ for every interval $I$ on the circle. Let $U_T$ be the Koopman operator. If there exists a spectral measure $\mu_f$ of $U_T$ which is uniformly $f$-continuous, then $T$ is partially weakly $f$-mixing. If $T$ is partially weakly $f$-mixing, then there exists $\mu_f$ of $U_T$ which is uniformly $\sqrt{f}$-continuous \cite{spec}. Carvalho and de Oliveira \cite{Carvalho} proved some properties of $\limsup_N N^{\alpha}C_N$ and $\liminf_N N^{\alpha}C_N$ for $0<\alpha<1$. Shortly after the first version of this paper was posted on arXiv, \cite{moll2023speed} proved quantitative weak-mixing estimates for the Chacon transformation for Lipschitz observables tested against $L^2$ observables. 

\subsection{Main Ideas and Outline of Paper}
\label{subsec:proof-sketch}
In Section~\ref{sec:exceptional-set}, we prove Proposition~\ref{prop:weak-conv-rate} and proceed to give a general method for constructing exceptional sets. Our main idea is the following. Given a measure-preserving transformation $(X, \mcr{B}, \mu, T)$ and $A \in \mcr{B}$, let $t_l: A \to \mathbb{N}_{\ge 0}$ denote the $l$-th return time of $x \in A$. That is,
\[ 
t_0(x) = 0, \quad t_1(x) = \min_{n \ge 1} \{ T^n(x) \in A \}, \quad t_{l+1}(x) = t_l(x) + t_1(T^{t_l}(x)).
\]
Now define $d_l$ to be the density of points with $l$-th return time $n$:
\[ 
d_l(n) = \mu \left( \{ x \in A: t_l(x) = n \}\right)= \mu(t_l^{-1} (n)).
\]
Note that $d_l$ has total mass $\mu(A)$ for each $l$. Clearly, this allows us to write
\[ 
\mu(A \cap T^{-n}A) = \mu(\{ x \in A : T^n(x) \in A \}) = \sum_{l = 0}^\infty d_l(n).
\]
If we normalize $d_l$ and extend it to $\mathbb{R}$ so that $D_l(x) = \mu(A)^{-1} d_l\left( \left\lfloor x + \frac{l}{\mu(A)} + \frac{1}{2} \right\rfloor \right)$, we have
\[ 
\sum_{l=0}^\infty d_l(n) = \mu(A) \sum_{l=0}^\infty D_l \left(n - \frac{l}{\mu(A)} \right).
\]
Now, if all of the $D_l$'s have roughly the same distribution $D_l \approx D$ (in a sense to be described later), then
\begin{equation}
\label{eq:heuristic}
\sum_{l=0}^\infty D_l \left(n - \frac{l}{\mu(A)} \right) \approx \sum_{l = 0}^\infty D \left(n - \frac{l}{\mu(A)} \right) \approx \mu(A) \int_{-\infty}^{\infty} D (x) \mathrm dx = \mu(A),
\end{equation}
which implies $\mu(A \cap T^{-n}(A)) \to \mu(A)^2$. Therefore, the exceptional sets for $(A, A)$ should contain the values of $n$ such that $\{ D_l : d_l(n) > 0 \}$ deviate significantly from $D$. Once we have exceptional sets of the form $J_{A, A}$, it is routine to generalize to other pairs.

In order to use this idea to actually construct exceptional sets, we need a way of describing $D_l$. In Section~\ref{sec:chacon}, we do this for the Chacon transformation using the recursive equation
\[ 
D_{3l}(x) = D_l(x) , \quad D_{3l \pm 1} = \frac{1}{3}\left\{ D_{l \pm 1}(x) + D_l \left( x + \frac{1}{2} \right) + D_l \left( x - \frac{1}{2} \right) \right\}.
\]
Note that the recursive equation for $D_{3l \pm 1}$ is a linear combination of $D_l, D_{l \pm 1}$ convolved by a probability measure. Furthermore, by repeating this recursion $q$ times, we see that $D_{3^q l}, D_{3^q l + 1}, \dots, D_{3^q(l+1)}$ are all convex combinations of $D_l, D_{l+1}$ and their translations by at most $q/2$ to each side. Thus, as long as $\| D_l - D_{l+1} \|_1$ and $\|D_l(x) - D_l(x-t)\|_1$ are small (for fixed $t$), we can argue along the lines of equation \eqref{eq:heuristic}. Because $D_l$ is symmetric and unimodal, both of these terms can be related to the peak value of $D_l$, which is often referred to as its height. If we let $b_l$ be the size of the support of $D_l$, we can notice the following facts. First, if $b_l$ is large, this implies that $D_l$ underwent many convolutions. Thus, by an application of the local limit theorem, we can deduce that $D_l(0)$ decays sufficiently. Secondly, if $b_l$ is small, it can only affect a small number of $n$ with $d_l (n) > 0$. Furthermore, the recursive equation for $b_l$ ensures that $b_l$ diverges to infinity except for a small number of values $l$. Thus, we can choose the values of $n$ in the support of $d_l$ with small $b_l$, and this will give us an exceptional set for $(A, A)$.

In Section~\ref{sec:rank-one}, we generalize to a larger class of transformations which we call \emph{restrictive tight maps}. Essentially, they are a class of transformations for which we can derive a recursive equation of the form
\[ 
D_{ml + r} = \frac{m-r}{m} \alpha_r \ast D_l + \frac{r}{m} \beta_r \ast D_{l+1}
\]
where for $1\le r<m$, $\alpha_r$, $\beta_r$ are probability distributions on $\frac{1}{m-1}\mathbb{Z}$. Clearly, the Chacon transformation is an example of a restrictive tight map\footnote{As such, Section~\ref{sec:chacon} is not strictly necessary. However, the Chacon case is much simpler and helps motivate many of the changes in the general setting, so we have decided to include it in our paper}. Given such an equation, we can proceed in a similar fashion as the Chacon map, but with some key changes. The first issue is that $D_l$ is no longer symmetric or unimodal. As such, knowing the height of $D_l$ is no longer sufficient to derive the necessary properties. To combat this, we turn to the \emph{total variation} of $D_l$, which controls the fluctuations of $D_l$ as well as its height. Of course, this means that the local central limit theorem is no longer sufficient, so we introduce a new proof using a coupling of random walks. The second problem has to do with $b_l$. In the Chacon case, $b_l$ gave us information on two important quantities: the size of the support of $D_l$, and how 'mixed' it is (i.e., the amount of convolutions applied to it). In general, this relationship is more delicate. As such, we introduce a new quantity $c_l$ that better describes the total variation, and then compare its size to $b_l$.

In Section~\ref{sec:main-proof}, we use the properties proven in previous sections to construct upper-bound exceptional sets for restrictive tight maps. The main counting tool is a quantitative estimate for numbers with small $c_l$. The same counting estimates also yield lower bounds for tight maps with no spacers above the last subcolumn, leading to the proofs of Theorems~\ref{thm:general} and~\ref{thm:lower-bound}. Corollary~\ref{cor:chacon} follows from the fact that the Chacon transformation is a restrictive tight map.

Finally, in Section~\ref{sec:applications}, we discuss generalizations and applications of our results. While weakly mixing systems are most commonly defined for discrete-time $\mathbb{Z}$-actions, these definitions can be generalized to continuous or higher-rank actions. We show that Proposition~\ref{prop:weak-conv-rate} can be easily adapted to these settings. We then apply these results, under the regularity hypotheses of the cited quantitative estimates, to random substitution tilings, interval exchange transformations, translation flows, primitive substitution $\mathbb{Z}$-actions, and self-affine substitution tilings. Subsection~\ref{subsec:rank-one-lower-bound} gives a separate weakly mixing one-spacer rank-one example for which exceptional sets for a suitable pair have polynomial lower growth. Lastly, we conclude with some related open questions in Section~\ref{subsec:open-questions}.

\subsection{Definitions and Notation}
\label{subsec:notation}
$(X, \mcr{B}, \mu, T)$ will denote a measure-preserving system. For any $A \in \mcr{B}$, let $a = \mu(A)^{-1}$. $r_A : A \to \mathbb{N}$ denotes the first return time $r_A(x) = \min \{ n \ge 1 : T^n(x) \in A \}$, and the returning point is given by $S_A(x) = T^{r_A(x)}(x)$. $t_l(x)$ refers to the $l$-th return time, i.e., $t_1 = r_A$ and $t_{l+1}(x) = t_1(x) + t_l (S_A(x))$. $d_l(n) = \mu(t_l^{-1}(n))$ has total mass $\mu(A)$ and $D_l(x) = a d_l (\lfloor x + la + 1/2 \rfloor)$. $C$, $C'$, $c$, and so on will refer to constants that may change from line to line. They may depend on the transformation $T$, but do not depend on the set $A \in \mcr{B}$, except in arguments where a particular $A_k$ has been fixed. 

\section{Weak Mixing and Exceptional Sets}
\label{sec:exceptional-set}

\subsection{\texorpdfstring{Proof of Proposition~\ref{prop:weak-conv-rate}}{Proof of Proposition}}
In this section, we prove Proposition~\ref{prop:weak-conv-rate} and show that the rate of weak mixing provides information on the size of exceptional sets. This is done through the following lemma, which is a modification of \cite[Lemma~2.41]{ward}.

\begin{lemma}
\label{lem:weak-conv-speed}
 Let $(a_n)$ be a bounded sequence of non-negative real numbers. Suppose that 
 \[
 \frac{1}{n} \sum_{j = 0}^{n-1} a_j = o(b_n)
 \]
 and $b_n \to 0$ as $n \to \infty$. Then, there exists a set $J \subseteq \mathbb{N}$ such that $|J \cap [0,n]| / n b_n$ converges to zero, and $a_n \to 0$ as $n \to \infty$ for all $n \notin J$.
\end{lemma}

\begin{proof}
 Define $J_k$ as
 $$ J_k = \left\{ j \in \mathbb{N} : a_j > \frac{1}{k} \right\} $$
 for all $k \in \mathbb{N}$. Then, it is clear that $J_1 \subseteq J_2 \subseteq \cdots $ and, for the half-open interval,
 $$\frac{1}{k} \left| J_k \cap [0,n) \right| \le \sum_{\substack{a_j > \frac{1}{k} \\ j < n}} a_j \le \sum_{j = 0}^{n-1} a_j = o( nb_n).$$
 The same estimate holds with $[0,n]$ in place of $[0,n)$. Indeed, if $J_k$ is empty this is clear; if $J_k$ is nonempty and finite, the hypothesis forces $(nb_n)^{-1}\to0$; and if $J_k$ is infinite, the half-open estimate above implies $(nb_n)^{-1}\to0$ after the first element of $J_k$. Hence,
 $$\left|J_k\cap[0,n]\right|=o(nb_n).$$
 Thus, we can define a strictly increasing sequence of positive integers $(l_k)$ such that
 $$ \frac{1}{nb_n} \left| J_k \cap [0,n] \right| \le \frac{1}{k} $$
 for all $n \ge l_k$ and $k \ge 1$. Now let us define $J$ as
 $$ J = \bigcup_{k = 1}^{\infty} \left( J_k \cap [l_k, l_{k+1}) \right) $$ 
 and show that $J$ satisfies our conditions.
 
 To see that $a_n \to 0$ as $n \to \infty$ and $n \notin J$, it is enough to notice that since $J_k \cap [l_k , \infty) \subseteq J$, $a_n \le \frac{1}{k}$ if $n \ge l_k$ and $n \notin J$. Further, since $J \cap [0,n] \subseteq J_k \cap [0,n]$ if $n \in [l_k, l_{k+1})$, $$ \frac{1}{nb_n} \left| J \cap [0,n] \right| \le \frac{1}{nb_n} \left| J_k \cap [0,n] \right| \le \frac{1}{k}$$ and so $\frac{1}{nb_n}|J \cap [0,n]|$ converges to zero.
\end{proof}

\begin{proof}[Proof of Proposition~\ref{prop:weak-conv-rate}]
Apply Lemma~\ref{lem:weak-conv-speed} to
\[
 a_n=\left|\mu(A\cap T^{-n}B)-\mu(A)\mu(B)\right|^p.
\]
The hypothesis gives \(N^{-1}\sum_{n<N}a_n=o(b_N)\), so there is
\(J_{A,B}\subseteq\mathbb N\) with
\(|J_{A,B}\cap[0,N]|=o(Nb_N)\) and \(a_n\to0\) outside \(J_{A,B}\). Hence
\(a_n^{1/p}\to0\) outside \(J_{A,B}\), which is the desired conclusion.
\end{proof}

Hence, we can find an upper bound on the size of the exceptional set given the rate of weak mixing.

\subsection{Constructing Exceptional Sets}

Now we provide a general strategy for constructing exceptional sets. This idea will be used in future sections to obtain exceptional sets for the Chacon transformation and restrictive tight maps. The following proposition is one of the key ideas of this paper.

\begin{proposition}
\label{prop:Fn-Gn}
Let \(A\in\mcr B\) satisfy \(\mu(A)>0\), put \(a=\mu(A)^{-1}\), and set \(D_l\equiv0\) for \(l<0\). For each \(n\), choose \(F_n,G_n\in BV(\mathbb R)\cap L^1(\mathbb R)\) such that
\[
0\le F_n\le G_n,\qquad \|F_n\|_1\le1\le\|G_n\|_1,
\]
and
\[
F_n(n-al)\le D_l(n-al)\le G_n(n-al)
\qquad(l\in\mathbb Z).
\]
Then, for any \(\epsilon(n)\to0\),
\[
\{ n\in\mathbb N:
\max(\|G_n-F_n\|_1,V(G_n),V(F_n))\ge\epsilon(n)\}
\]
is an exceptional set for \((A,A)\).
\end{proposition}

\begin{proof}
For integer \(n\), the normalization gives \(D_l(n-al)=a d_l(n)\), and hence
\[
\mu(A\cap T^{-n}A)=\sum_{l\ge0}d_l(n)=\frac1a\sum_{l\in\mathbb Z}D_l(n-al).
\]
Using the pointwise bounds and the Riemann-sum estimate from the Appendix, we obtain
\[
\mu(A)^2\|F_n\|_1-\mu(A)V(F_n)
\le
\mu(A\cap T^{-n}A)
\le
\mu(A)^2\|G_n\|_1+\mu(A)V(G_n).
\]
Outside the displayed exceptional set, the two variations and \(\|G_n-F_n\|_1\) tend to zero. Since
\(\|F_n\|_1\le1\le\|G_n\|_1\) and
\(\|G_n\|_1-\|F_n\|_1\le\|G_n-F_n\|_1\), both norms tend to \(1\). The sandwich estimate then gives
\(\mu(A\cap T^{-n}A)\to\mu(A)^2\).
\end{proof}

\begin{remark}
 A natural choice of $F_n, G_n$ would be $F_n = \min_{l \in P_n} D_l$ and $G_n = \max_{l \in P_n} D_l$. However, it is not always easy to show that these maps have small total variation. In future sections, we will choose appropriate $F_n$ and $G_n$ as necessary.
\end{remark}

Now we explain how we can construct exceptional sets using preexisting ones. We say that a $J_1 \subseteq \mathbb{N}$ \emph{eventually contains} $J_2 \subseteq \mathbb{N}$ if $J_2 \setminus J_1$ is finite, i.e., $J_1$ contains all but finitely many elements of $J_2$. It is clear that if $J_{A, B}$ is an exceptional set, then any set that eventually contains $J_{A, B}$ is also exceptional. We state this fact in the following lemma for future reference.

\begin{lemma}\label{lem:finite-modification-exceptional}
Let $(X,\mathscr B,\mu,T)$ be a measure-preserving system. Let $A,B\in\mathscr B$ and $J \subseteq \mathbb N$. 
\begin{enumerate}
\item If $J_{A,B}\setminus J$ is finite for some exceptional set $J_{A,B}$, then $J$ is exceptional for $(A,B)$.

\item Let $A\in\mathscr B$, and for $\tau>0$ define
\[
\mathcal E_A^T(\tau):=\{n\in\mathbb N:\ |\mu(A\cap T^{-n}A)-\mu(A)^2|>\tau\}.
\]
Then $\mathcal E_A^T(\tau)\setminus J_{A,A}$ is finite.
\end{enumerate}
\end{lemma}

\begin{proof}
Both assertions are immediate from the definitions.
\end{proof}

Now by the following lemma, we may generate exceptional sets for any countable collection of measurable sets.

\begin{lemma}
\label{lem:eventual-containment}
Suppose $J_1,J_2,\dots\subseteq\mathbb N$ such that, for each $i$, the estimate
\[
|J_i\cap[0,n]|\le f(n)
\]
holds for all sufficiently large $n$. Then for any increasing function $h(n)$ diverging to infinity, we can construct a set $J$ such that
\[
|J\cap[0,n]|\le h(n)f(n)
\]
for all sufficiently large $n$, and every $J_i$ is eventually contained in $J$.
\end{lemma}

\begin{proof}
For each $i$, choose $N_i$ so that $|J_i\cap[0,n]|\le f(n)$ for all $n\ge N_i$, and set $J_i'=J_i\cap[N_i,\infty)$. Then $|J_i'\cap[0,n]|\le f(n)$ for every $n$, and replacing $J_i$ by $J_i'$ does not change eventual containment. Applying the construction in Theorem~3.1 of \cite{Friedman} to the sequence $J_i'$ gives
\[
J=\bigcup_{i=1}^{\infty}\left(J_i'\setminus\{m\in\mathbb N:h(m)\le i\}\right).
\]
Then $|J\cap[0,n]|\le h(n)f(n)$ for all sufficiently large $n$, and every $J_i'$ is eventually contained in $J$. Hence every original $J_i$ is eventually contained in $J$.
\end{proof}

\begin{corollary}
\label{prop:J_C}
 Let $\mcr{C} \subseteq \mcr{B}$ be a countable collection of measurable sets. Suppose that for any $A, B \in \mcr{C}$, there exists an exceptional set $J_{A, B}$ of $(A, B)$. Further, assume that, for each $A, B \in \mcr{C}$, the estimate $|J_{A, B} \cap [0, n]| \le f(n)$ holds for all sufficiently large $n$. Then, given any increasing function $h(n)$ diverging to infinity, we can construct an exceptional set $J_{\mcr{C}}$ of $\mcr{C}$ such that $|J_{\mcr{C}} \cap [0, n]| \le f(n) h(n)$ for all sufficiently large $n$. 
\end{corollary}

\begin{proof}
    The proof is immediate from Lemma~\ref{lem:eventual-containment}.
\end{proof}

The next proposition shows that an exceptional set for a generator is also exceptional for the entire $\sigma$-algebra.
\begin{proposition}
 \label{prop:J_B}
Let $\mcr{C} \subseteq \mcr{B}$ be a countable generator of $\mcr{B}$ and suppose $J_{\mcr{C}}$ is exceptional for $\mcr{C}$. Further, suppose for every $\varepsilon>0$ and $A \in \mcr{B}$, there exist finitely many mutually disjoint sets $A_i \in \mcr{C}$ such that 
 \begin{equation*}
 \mu(A\Delta (\cup_i A_i))<\varepsilon.
 \end{equation*}
 Then, $J_{\mcr{C}}$ is exceptional for $\mcr{B}$.
\end{proposition}

\begin{proof}
  The proof is a standard approximation via generating sets and goes along the same lines as Theorem~5.11 of \cite{Friedman2}. Hence we omit the proof here, and refer the interested reader to \cite{Friedman2}.
\end{proof}

In most of this paper we define exceptional sets in the setting of measurable sets, that is, by studying
\[
\mu\bigl(A\cap T^{-n}B\bigr).
\]
However, one can extend this notion to general $L^2$ functions as follows. Let \(f,g\in L^2(\mu)\) and write \(\mu(f)=\int f\,d\mu\). Notice that
\begin{enumerate}
 \item \(\mu(A)=\int\chi_A\,d\mu=\mu(\chi_A)\),
 \item \(\mu(A\cap B)=\int\chi_A\,\chi_B\,d\mu=\mu(\chi_A\,\chi_B)\),
 \item \(\mu\bigl(T^{-n}A\bigr)=\int\chi_{T^{-n}A}\,d\mu
 =\int(\chi_A\circ T^n)\,d\mu
 =\mu(\chi_A\circ T^n).
\)
\end{enumerate}
Hence the correlation
\[
\mu\bigl(\chi_A\cdot(\chi_B\circ T^n)\bigr)
\;=\;\mu\bigl(A\cap T^{-n}B\bigr)
\]
and more generally one may consider
\[
\mu\bigl(f\cdot (g\circ T^n)\bigr)
\;=\;\int f\,(g\circ T^n)\,d\mu.
\]
We then say a set \(J_{f,g}\subseteq\mathbb{N}\) is \emph{exceptional} for the pair \((f,g)\) if
\[
\mu\bigl(f\cdot(g\circ T^n)\bigr)\to\mu(f)\,\mu(g)
\quad(n\to\infty,\;n\notin J_{f,g}).
\]
In particular, when \(f=\chi_A\) and \(g=\chi_B\) this recovers the usual definition, since
\(\chi_A\cdot(\chi_B\circ T^n)=\chi_{A\cap T^{-n}B}\), and thus \(J_{\chi_A,\chi_B}=J_{A,B}\).

\begin{proposition}
 \label{prop:J_f}
 Let $f, g\in L^2(\mu)$. If $J$ is exceptional for $\mcr{B}$, then $J$ is exceptional for $(f, g)$.
\end{proposition}

\begin{proof}
Again, the proof is a standard approximation via simple functions. Let 
\[
f_m=\sum_{i=1}^{k_m}\alpha_i^{(m)}\chi_{A_i^{(m)}},
\qquad
g_m=\sum_{j=1}^{\ell_m}\beta_j^{(m)}\chi_{B_j^{(m)}}
\]
be simple functions with 
\(\|f-f_m\|_{2}<\delta_m\) and \(\|g-g_m\|_{2}<\delta_m\), where \(\delta_m\to0\). 
Since \(T\) is measure preserving,
\(\|g\circ T^n\|_{2}=\|g\|_{2}\) and 
\(\|(g_m-g)\circ T^n\|_{2}=\|g_m-g\|_{2}\) for all \(n\). 

Fix \(\epsilon>0\) and choose \(m\) sufficiently large so that
\[
|\mu(f_m)\,\mu(g_m)-\mu(f)\,\mu(g)|<\frac{\epsilon}{3},
\quad
\|f-f_m\|_{2}\,\|g\|_{2}+\|f_m\|_{2}\,\|g_m-g\|_{2}<\frac{\epsilon}{3}.
\]
For simple functions, the correlation error is a finite linear combination of indicator-pair correlation errors, so \(J\) is exceptional for the simple pair \((f_m,g_m)\). Hence there is \(N\) such that for all \(n>N\), \(n\notin J\),
\[
\bigl|\mu\bigl(f_m\,(g_m\circ T^n)\bigr)-\mu(f_m)\,\mu(g_m)\bigr|<\frac{\epsilon}{3}.
\]
Moreover for any such \(n\),
\[
\bigl|\mu\bigl(f\,(g\circ T^n)\bigr)-\mu\bigl(f_m\,(g_m\circ T^n)\bigr)\bigr|
\le
\|f-f_m\|_{2}\,\|g\|_{2}
+\|f_m\|_{2}\,\|g_m-g\|_{2}
<\frac{\epsilon}{3}.
\]
Hence for all \(n>N\), \(n\notin J\),
\begin{align*}
\bigl|\mu\bigl(f\,(g\circ T^n)\bigr)-\mu(f)\,\mu(g)\bigr|
&\le
{\bigl|\mu(f(g\circ T^n))-\mu(f_m(g_m\circ T^n))\bigr|}\\
&\quad+
{\bigl|\mu(f_m(g_m\circ T^n))-\mu(f_m)\,\mu(g_m)\bigr|}\\
&\quad+
{|\mu(f_m)\,\mu(g_m)-\mu(f)\,\mu(g)|}\\
&<\epsilon.
\end{align*}
Therefore \(J\) is exceptional for \((f,g)\).
\end{proof}

\section{The Chacon Transformation}
\label{sec:chacon}
\subsection{The Chacon Transformation}
\label{sec:chacon-def}

\begin{figure}[b!]
\centering
 \begin{tikzpicture} [xscale = 5, yscale = 5,circ/.style={draw,circle,inner sep=0pt,minimum size=4pt},
 dot/.style={circ,fill}]
 
 \node[black] at (-0.6,0) {\Large $\text{Step } 0$};
 \node[black] at (0.3333, 0.3) {\Large$\text{tower}$};
 \node[black] at (1.2666, 0.3) {\Large $\text{spacer}$};
 \draw[black, ultra thick] (0,0) -- (0.6666,0);
 \draw[blue, ultra thick] (1.1,0) -- (1.4333, 0);
 \node[black] at (-0.06, 0) {$0$};
 \node[black] at (0.7606, 0) {$2/3$};
 \node[blue] at (1.0, 0) {$2/3$};
 \node[blue] at (1.4933, 0) {$1$};
 \draw[thick, dotted, black] (0.2222, -0.05)--(0.2222, 0.05);
 \draw[thick, dotted, black] (0.4444, -0.05)--(0.4444, 0.05);
 \draw[thick, dotted, black] (1.3222, -0.05)--(1.3222, 0.05);
 \draw[gray, densely dashed] (0.1111, 0.03) --(0.1111, 0.10)--(-0.2, 0.10)--(-0.2, -0.15)--(0.3333, -0.15);
 \draw[gray, -stealth, densely dashed] (0.3333, -0.15)--(0.3333, -0.03);
 \draw[gray, densely dashed] (0.3333, 0.03)--(0.3333, 0.10)--(1.2111, 0.10);
 \draw[gray, -stealth, densely dashed] (1.2111, 0.10)--(1.2111, 0.03);
 \draw[gray, densely dashed] (1.2111, -0.03) --(1.2111, -0.15)--(0.5555, -0.15);
 \draw[gray, -stealth, densely dashed] (0.5555, -0.15)--(0.5555, -0.03);
 
 \node[black] at (-0.6,-0.75) {\Large $\text{Step } 1$};
 
 \draw [blue, ultra thick] (0, -0.6)--(0.6666, -0.6);
 \draw [black, ultra thick] (0, -0.45)--(0.6666, -0.45);
 \draw [black, ultra thick] (0, -0.75)--(0.6666, -0.75);
 \draw [black, ultra thick] (0, -0.9)--(0.6666, -0.9);
 
 \draw[blue, ultra thick] (1.1, -0.75)--(1.4333, -0.75);
 
 \draw[thick, dotted, black] (0.2222, -0.95)--(0.2222, -0.4);
 \draw[thick, dotted, black] (0.4444, -0.95)--(0.4444, -0.4);
 \draw[thick, dotted, black] (1.3222, -0.8)--(1.3222, -0.7);
 
 \draw[black, -stealth] (0.1111, -0.87)--(0.1111, -0.78);
 \draw[black, -stealth] (0.3333, -0.87)--(0.3333, -0.78);
 \draw[black, -stealth] (0.5555, -0.87)--(0.5555, -0.78);
 \draw[black, -stealth] (0.1111, -0.72)--(0.1111, -0.63);
 \draw[black, -stealth] (0.3333, -0.72)--(0.3333, -0.63);
 \draw[black, -stealth] (0.5555, -0.72)--(0.5555, -0.63);
 \draw[black, -stealth] (0.1111, -0.57)--(0.1111, -0.48);
 \draw[black, -stealth] (0.3333, -0.57)--(0.3333, -0.48);
 \draw[black, -stealth] (0.5555, -0.57)--(0.5555, -0.48);
 
 \draw[gray, densely dashed] (0.1111, -0.42)--(0.1111, -0.35)--(-0.2, -0.35)-- (-0.2, -1.05)--(0.3333, -1.05);
 \draw[gray, -stealth, densely dashed] (0.3333, -1.05)--(0.3333, -0.93);
 
 \draw[gray, densely dashed] (0.3333, -0.42)--(0.3333, -0.35)--(1.2111, -0.35);
 \draw[gray, -stealth, densely dashed] (1.2111, -0.35)--(1.2111, -0.72);
 \draw[gray, densely dashed] (1.2111, -0.78)--(1.2111, -1.05)--(0.5555, -1.05);
 \draw[gray, -stealth, densely dashed] (0.5555, -1.05)--(0.5555, -0.93);
 
 \node[black] at (-0.06, -0.9) {$0$};
 \node[black] at (0.7666, -0.9) {$2/9$};
 \node[black] at (-0.1, -0.75) {$2/9$};
 \node[black] at (0.7666, -0.75) {$4/9$};
 \node[blue] at (-0.1, -0.6) {$2/3$};
 \node[blue] at (0.7666, -0.6) {$8/9$};
 \node[black] at (-0.1, -0.45) {$4/9$};
 \node[black] at (0.7666, -0.45) {$2/3$};
 
 \node[blue] at (1.0, -0.75) {$8/9$};
 \node[blue] at (1.4933, -0.75) {$1$};
 
 \draw[black,thick, fill=black] (0,0) circle (0.4pt);
 \draw[black,thick, fill=white] (0.6666,0) circle (0.4pt);
 \draw[black,thick, fill=black] (0,-0.45) circle (0.4pt);
 \draw[black,thick, fill=white] (0.6666,-0.45) circle (0.4pt);
 \draw[blue,thick, fill=blue] (0,-0.6) circle (0.4pt);
 \draw[blue,thick, fill=white] (0.6666,-0.6) circle (0.4pt);
 \draw[black,thick, fill=black] (0,-0.75) circle (0.4pt);
 \draw[black,thick, fill=white] (0.6666,-0.75) circle (0.4pt);
 \draw[black,thick, fill=black] (0,-0.9) circle (0.4pt);
 \draw[black,thick, fill=white] (0.6666,-0.9) circle (0.4pt);
 \draw[blue,thick, fill=blue] (1.1,0) circle (0.4pt);
 \draw[blue,thick, fill=white] (1.4333,0) circle (0.4pt);
 \draw[blue,thick, fill=blue] (1.1,-0.75) circle (0.4pt);
 \draw[blue,thick, fill=white] (1.4333,-0.75) circle (0.4pt);
 
\end{tikzpicture} 
\caption{Construction of the Chacon transformation in step $0$ and $1$ ($\tau_1 $)}
\label{Chacon_fig}
 \end{figure}
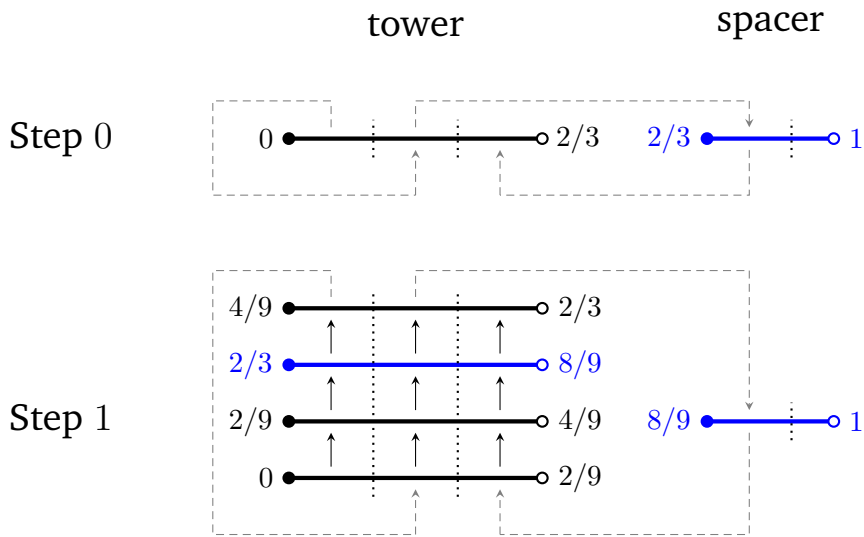 

In this section, we define the Chacon transformation. We start with two intervals, $[0,\frac{2}{3})$ and $[\frac{2}{3}, 1)$. The first interval is called the \emph{tower}, while the second is called the \emph{spacer}. (See Figure~\ref{Chacon_fig}) Then, we cut the tower into three pieces of equal width, $[0,\frac{2}{9}), [\frac{2}{9}, \frac{4}{9}),$ and $[\frac{4}{9}, \frac{2}{3})$. We also cut the spacer into two pieces, $[\frac{2}{3}, \frac{8}{9})$ and $[\frac{8}{9}, 1)$. Note that the width of the first piece of the spacer is double the width of the second, and is also equal to the width of each piece of the tower. Now we ``stack'' these pieces so that the first piece of the tower goes on the bottom, then the second piece, then the first piece of the spacer, and then the third piece of the tower (See Figure~\ref{Chacon_fig}). $\tau_1 : [0, \frac{4}{9}) \cup [\frac{2}{3}, \frac{8}{9}) \to [\frac{2}{9}, \frac{8}{9})$ is the map that sends each point in the tower (besides those on the top) to the point directly above it. This map is represented by the solid arrows in Figure~\ref{Chacon_fig}. For instance, $\tau_1 (\frac{1}{3}) = \frac{7}{9}$.

After this first step, we have a tower of height $4$ and width $\frac{2}{9}$ and a spacer of width $\frac{1}{9}$. Now we repeat this process. At every step, we cut the tower into three equal pieces, and the spacer into two, so that the width of the first piece of the spacer matches the width of the pieces cut from the tower. Then we stack the intervals in the same order as before. More precisely, we stack the middle third on top of the left third, then the first piece of the spacer, and lastly place the right third on top (See Figure~\ref{Chacon_fig2}). This implies that the spacer will be placed roughly one third of the way from the top. Then, we can observe that after the $n$th step:

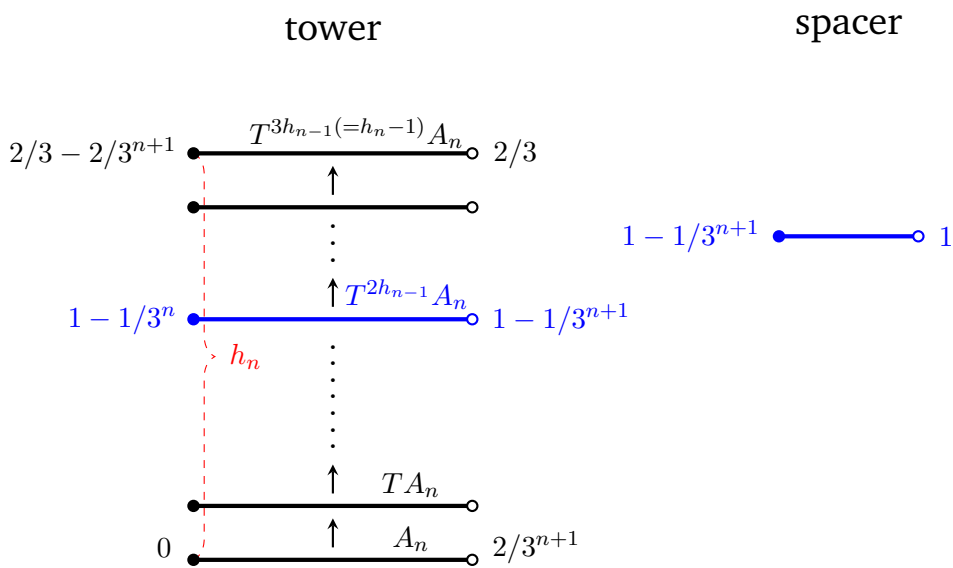
\begin{figure} [h]
\centering
 \begin{tikzpicture} [xscale = 5.5, yscale = 5.5,
 circ/.style={draw,circle,inner sep=0pt,minimum size=4pt},
 dot/.style={circ,fill}
 ]
\draw [red, dashed, decorate,decoration={brace,amplitude=8pt}]
 (0,0) -- (0,-0.98) node [midway, xshift=0.27in] {$h_n$};
 
 \node[black] at (0.3333, 0.3) {\Large$\text{tower}$};
 \node[black] at (1.5666, 0.3) {\Large $\text{spacer}$};
 
 \draw[ultra thick, black] (0,0) --(0.6666, 0);
 \draw[ultra thick, black] (0, -0.13) -- (0.6666, -0.13);
 \draw[ultra thick, blue] (0, -0.40) -- (0.6666, -0.40);
 \draw[ultra thick, black] (0, -0.85)--(0.6666, -0.85);
 \draw[ultra thick, black] (0, -0.98)--(0.6666, -0.98);

 \draw[black, thick, -stealth] (0.3333, -0.1) -- (0.3333, -0.03);
 \draw[black, thick, -stealth] (0.3333, -0.37) -- (0.3333, -0.30);
 \draw[black, thick, -stealth] (0.3333, -0.95) -- (0.3333, -0.88);
 \draw[black, thick, -stealth] (0.3333, -0.82) -- (0.3333, -0.75);
 
 \node [black, scale=0.45] at (0.3333,-0.26) {\textbullet};
 \node [black, scale=0.45] at (0.3333,-0.22) {\textbullet};
 \node [black, scale=0.45] at (0.3333,-0.18) {\textbullet};

 \node [black, scale=0.45] at (0.3333,-0.71) {\textbullet};
 \node [black, scale=0.45] at (0.3333,-0.67) {\textbullet};
 \node [black, scale=0.45] at (0.3333,-0.63) {\textbullet};
 \node [black, scale=0.45] at (0.3333,-0.59) {\textbullet};
 \node [black, scale=0.45] at (0.3333,-0.55) {\textbullet}; 
 \node [black, scale=0.45] at (0.3333,-0.51) {\textbullet};
 \node [black, scale=0.45] at (0.3333,-0.47) {\textbullet};
 
 \node[black] at (-0.24, 0) {$2/3-2/3^{n+1}$};
 \node[black] at (0.77, 0) {$2/3$};
 \node[black] at (-0.07, -0.95) {$0$};
 \node[black] at (0.82, -0.95) {$2/3^{n+1}$};
 
\node[black] at (0.52, -0.93) {$ A_n$};
\node[black] at (0.52, -0.8) {$TA_n$};
\node[blue] at (0.51, -0.34) {$ T^{2h_{n-1}}A_n$};

\node[black] at (0.39, 0.05) {$T^{3h_{n-1}(=h_n-1) }A_n$};
 
 \node[blue] at (-0.17, -0.4) {$1-1/3^{n}$};
 \node[blue] at (0.88, -0.4) {$1-1/3^{n+1}$};
 
 \draw[blue, ultra thick] (1.4, -0.2)-- (1.7333, -0.2);
 \node[blue] at (1.19, -0.2) {$1-1/3^{n+1}$};
 \node[blue] at (1.8, -0.2) {$1$};
 
\draw[black,thick, fill=black] (0,0) circle (0.35pt);
 \draw[black,thick, fill=white] (0.6666,0) circle (0.35pt);
 \draw[black,thick, fill=black] (0,-0.13) circle (0.35pt);
 \draw[black,thick, fill=white] (0.6666,-0.13) circle (0.35pt);
 \draw[blue,thick, fill=blue] (0,-0.4) circle (0.35pt);
 \draw[blue,thick, fill=white] (0.6666,-0.4) circle (0.35pt);
 \draw[black,thick, fill=black] (0,-0.85) circle (0.35pt);
 \draw[black,thick, fill=white] (0.6666,-0.85) circle (0.35pt);
 \draw[black,thick, fill=black] (0,-0.98) circle (0.35pt);
 \draw[black,thick, fill=white] (0.6666,-0.98) circle (0.35pt);
 \draw[blue,thick, fill=blue] (1.4,-0.2) circle (0.35pt);
 \draw[blue,thick, fill=white] (1.7333,-0.2) circle (0.35pt);
 
\end{tikzpicture} 
\caption{Construction of the Chacon transformation in step $n\geq 2$ }
\label{Chacon_fig2}
 \end{figure} 

\begin{enumerate}
 \item The height of tower is $h_n=3h_{n-1}+1$, where $h_0 = 1$.
 \item The width of the spacer is $3^{-(n+1)}$.
 \item The width of each interval is $2\cdot 3^{-(n+1)}$.
\end{enumerate}

Note that $h_n=(3^{n+1}-1)/2$. As before, we define $$\tau_n : \left[0, 1\right) \setminus \bkt{\left[\frac{2}{3}- \frac{2}{3^{n+1}}, \frac{2}{3}\right) \cup \left[1 - \frac{1}{3^{n+1}}, 1\right)} \to \left[\frac{2}{3^{n+1}} , 1 - \frac{1}{3^{n+1}}\right)$$ to be the map sending each point in the tower to the one above it. Because of the way the tower is constructed, the values of the maps $\tau_n$ coincide whenever two of them are defined at the same point. Since the complements of the domains have summable measures $O(3^{-n})$, Borel--Cantelli, together with the countable set of endpoints, implies that $\tau_n(x)$ is eventually defined and the compatible values stabilize for a.e. $x$.

\begin{definition}
 The Chacon transformation $T:[0,1)\to [0,1)$ is defined by $T(x)=\lim_{n\to \infty} \tau_n(x)$ for Lebesgue-a.e. $x$, and is regarded modulo null sets.
\end{definition}

\begin{proposition}[\cite{Friedman,Alpha}]
 Let $T$ be the Chacon transformation defined above.
 \begin{enumerate}
 \item $T$ is measure preserving, ergodic, and weak mixing (with respect to the Lebesgue measure).
 \item $T$ is not mixing.
 \end{enumerate}
 \label{prop:chacon-prop}
\end{proposition}

Let $A_k := [0, 2/3^{k+1})$ be the bottom interval in the $k$-th step of the cutting and stacking process. Observing the cutting and stacking operations used to define the Chacon transformation, we see that, when ignoring the width and height of the towers, the same stacking procedure is applied at each step. As such, it is often useful to identify $A_k$ with the interval $[0,1)$ via the bijection $u_k : A_k \to [0,1)$ defined by $u_k(x) = a_kx$, where $a_k = \mu(A_k)^{-1}$. This allows us to disregard the change in width at every step. We use this bijection to redefine all the functions so that they are defined on $[0,1)$ rather than $A_k$:
\begin{gather}
 \label{eq:r_k}
 r_k : [0,1) \to \mathbb{N} \, \quad r_k = r_{A_k} \circ u_k ^{-1}, \\
 \label{eq:S_k}
 S_k : [0,1) \to [0,1) \, \quad S_k = u_k \circ S_{A_k} \circ u_k^{-1}, \\
 \label{eq:t_l'}
 t_l' : [0,1) \to \mathbb{N} \, \quad t_l' = t_l \circ u_k ^{-1}, \\
 \label{eq:d_l'}
 d_l '=d_l'(k) : \mathbb{N} \to \mathbb{R} \, \quad d_l'(n) =\mu \left( (t_l ')^{-1} (n) \right).
\end{gather}

Among these functions, $S_k$ is the only map that is independent of $k$ (see Lemma~\ref{lem:r_k}). However, as $k$ remains constant throughout all sections where the above functions are used, we have chosen to omit $k$ in our notation and use $r, S, t_l', d_l, d_l' \dots$ to denote the functions above. We can also see that $t_l ' (x) = \sum_{i = 0}^{l-1} r\left(S^{i}(x)\right)$ and $d_l ' = \frac{3^{k+1}}{2} d_l$ hold. Note that while $P_n$, $B_l$, and $b_l$ will be defined below with respect to $d_l$, we might as well have defined them using $d_l '$, as the two functions are simply scalar multiples of each other.

Lastly, because powers of $3$ come up often in our analysis, it is often convenient to write values using the ternary number system. From now on, we denote numbers using the ternary system. For instance, we have $\ovl{0.{2}} = \frac{2}{3}$ and $\ovl{0.{12}}=\frac{5}{9}$.

The main results of this section are Corollary~\ref{cor:D_l} and Lemma~\ref{lem:ternary}. Corollary~\ref{cor:D_l} gives a recurrence formula for $D_l$, which is the backbone of all future results. Further, Lemma~\ref{lem:ternary} shows that $b_l$, the size of the support of $D_l$, is related to the balanced ternary expansion of $l$ (see Lemma~\ref{lem:ternary}). Some of the results presented here have been shown in previous works such as \cite{Alpha}. Namely, parts of Corollary~\ref{cor:D_l} were shown in Theorem~3.3 and Proposition~4.2 in \cite{Alpha}.

\subsection{The Recursive Formula}
In this section, we derive a recursive formula for $D_l$ in the case of the Chacon transformation. Some of the results presented here have been shown in previous
 works such as \cite{Alpha}. Namely, Corollary~\ref{cor:D_l} is similar to \cite[Theorem~3.3]{Alpha}
\label{subsec:chacon-recursion}
\begin{lemma}
 \label{lem:r_k}
Let $r_k$ and $S = S_k$ be as in \eqref{eq:r_k} and \eqref{eq:S_k}. We have
 $$ r_k(\ovl{0.a_1 a_2 a_3 \cdots})=
 \begin{cases}
 h_k & \textrm{if } a_1 = 0 \\
 h_k + 1 & \textrm{if } a_1 = 1 \\
 r_k(\ovl{0.a_2 a_3 \cdots}) & \textrm{if } a_1 = 2,
 \end{cases}$$
 $$ S(\ovl{0.a_1 a_2 a_3 \cdots})=
 \begin{cases}
 \ovl{0.1 a_2 a_3 \cdots} & \textrm{if } a_1 = 0 \\
 \ovl{0.2 a_2 a_3 \cdots} & \textrm{if } a_1 = 1 \\
 \frac{1}{3}S(\ovl{0.a_2 a_3 \cdots}) & \textrm{if } a_1 = 2.
 \end{cases}$$
\end{lemma}

\begin{proof}
 Let us consider the position of $T^{h_k - 1} (\ovl{0.a_1 a_2 a_3 \cdots})$ with respect to the tower in step $k$.
 
 If $a_1 = 0$, then $T^{h_k - 1} (\ovl{0.a_1 a_2 a_3 \cdots})$ is at the left third of the topmost segment of the tower. Therefore, since the middle third of $A_k$ gets stacked above it in step $k+1$, $ r_k(\ovl{0.a_1 a_2 a_3 \cdots})= h_k$ and $ S(\ovl{0.a_1 a_2 a_3 \cdots})= \ovl{0.1 a_2 a_3 \cdots}$.
 
 Similarly, if $a_1 = 1$, then $T^{h_k - 1} (\ovl{0.a_1 a_2 a_3 \cdots})$ is at the middle third of the topmost segment of the tower. Therefore, since the spacer $A_k$ gets stacked above it, and the right third of $A_k$ above the spacer, $ r_k(\ovl{0.a_1 a_2 a_3 \cdots})= h_k + 1$ and $ S(\ovl{0.a_1 a_2 a_3 \cdots})= \ovl{0.2 a_2 a_3 \cdots}$.
 
 Lastly, consider the case where $a_1 = 2$. Then, $T^{h_k - 1} (\ovl{0.a_1 a_2 a_3 \cdots})$ is at the right third of the topmost segment of the tower. Thus, after Step $k+1$, it is still at the topmost segment of the tower. Furthermore, its position relative to the segment length is precisely $\ovl{0.a_2 a_3 \cdots}$. We also see that the segment length after step $k+1$ is one third of that at step $k$. Hence, we see that $ r_k(\ovl{0.a_1 a_2 a_3 \cdots})= r_k(\ovl{0.a_2 a_3 \cdots})$ and $ S(\overline{0.a_1 a_2 a_3 \cdots})= \frac{1}{3} S(\overline{0.a_2 a_3 \cdots})$.
\end{proof}

\begin{corollary}
 $$ S^{3l}(\ovl{0.a_1 a_2 a_3 \cdots})= \frac{1}{3} \left\{ a_1 + S^l ( \ovl{0. a_2 a_3 \cdots}) \right\} $$
 \label{cor:S}
\end{corollary}

\begin{proof}
 The cases where $l = 0, 1$ are either trivial or follow directly from Lemma~\ref{lem:r_k}, and all other cases follow naturally. (Note that $\frac{1}{3} (a_1 + \ovl{0.b_1 b_2 \cdots} ) = \ovl{0.a_1 b_1 b_2 \cdots}$.)
\end{proof}

\begin{lemma}
 \label{lem:t_l}
 Let $t_l'$ be as in \eqref{eq:t_l'}. Clearly, $t_0' = 0$. We have
 $$ t_{3l}'(\ovl{0.a_1 a_2 a_3 \cdots})= 2lh_k + l + t_l'(\ovl{0.a_2 a_3 \cdots}),$$
 
 $$ t_{3l + 1}'(\ovl{0.a_1 a_2 a_3 \cdots})=
 \begin{cases}
 (2l+1)h_k + l + t_l'(\ovl{0.a_2 a_3 \cdots}) & \textrm{if } a_1 = 0 \\
 (2l+1)h_k + l + 1 + t_l'(\ovl{0.a_2 a_3 \cdots}) & \textrm{if } a_1 = 1 \\
 2lh_k + l + t_{l+1}'(\ovl{0.a_2 a_3 \cdots}) & \textrm{if } a_1 = 2,
 \end{cases}$$
 
 $$ t_{3l + 2}'(\ovl{0.a_1 a_2 a_3 \cdots})=
 \begin{cases}
 (2l+2)h_k + l + 1 + t_l'(\ovl{0.a_2 a_3 \cdots}) & \textrm{if } a_1 = 0 \\
 (2l+1)h_k + l + 1 + t_{l+1}'(\ovl{0.a_2 a_3 \cdots}) & \textrm{if } a_1 = 1 \\
 (2l+1)h_k + l + t_{l+1}'(\ovl{0.a_2 a_3 \cdots}) & \textrm{if } a_1 = 2.
 \end{cases}$$
\end{lemma}

\begin{proof}
 
 Since all cases can be shown similarly, we only prove the first statement, $$ t_{3l}'(\ovl{0.a_1 a_2 a_3 \cdots})= 2lh_k + l + t_l'(\ovl{0.a_2 a_3 \cdots}).$$ In particular, let us focus on the case where $a_1 = 0$. Then, 
 \begin{equation*}
 \begin{split}
 t_{3l}' (\ovl{0.0 a_2 a_3 \cdots} ) &= \sum_{i = 0} ^{3l - 1} r( S^{i}(\ovl{0.0 a_2 a_3 \cdots}))\\
 &= \sum_{i = 0} ^{l-1} \left( r( S^{3i}(\ovl{0.0 a_2 a_3 \cdots})) + r( S^{3i + 1}(\ovl{0.0 a_2 a_3 \cdots})) + r( S^{3i + 2}(\ovl{0.0 a_2 a_3 \cdots})) \right) \\
 &= \sum_{i = 0} ^{l - 1} \left( r(\frac{1}{3} S^{i}(\ovl{0.a_2 a_3 \cdots})) + r(\frac{1}{3}(1 + S^{i}(\ovl{0.a_2 a_3 \cdots}))) + r(\frac{1}{3} (2 + S^{i}(\ovl{0.a_2 a_3 \cdots}))) \right) \\
 &= \sum_{i = 0} ^{l-1} \left( h_k + (h_k + 1) + r(S^i (\ovl{0.a_2 a_3 \cdots})) \right) \\
 &= 2lh_k + l + t_l'(\ovl{0.a_2 a_3 \cdots}).
 \end{split}
 \end{equation*}
 The third equality comes from Corollary~\ref{cor:S}, and the fourth comes from Lemma~\ref{lem:r_k}.
\end{proof}

\begin{lemma}
 \label{lem:d_l'}
 Let $d_l'$ be as in \eqref{eq:d_l'}. Clearly, $d_0' = \mathbf{1}_{\{ 0 \}}$. We have
 $$ d_{3l}'(i) = d_l' (i - 2lh_k - l)$$
 $$ d_{3l + 1}'(i)=
 \frac{1}{3} \left( d_l' (i - (2l + 1)h_k - l) + d_l'(i - (2l + 1)h_k - l - 1) + d_{l+1}'(i - 2lh_k - l) \right)$$
 $$ d_{3l + 2}'(i)=
 \frac{1}{3} \left( d_l' (i - (2l + 2)h_k - l - 1) + d_{l+1}'(i - (2l + 1)h_k - l - 1) + d_{l+1}'(i - (2l + 1)h_k - l) \right).$$
\end{lemma}

\begin{proof}
 Due to Lemma~\ref{lem:t_l}, $t_{3l}' (\ovl{0.a_1 a_2 \cdots}) = i$ if and only if $t_{l}' (\ovl{0.a_2 a_3 \cdots}) = i - 2lh_k -l$. Hence, $ d_{3l}'(i) = d_l' (i - 2lh_k - l)$. Similarly, $t_{3l+1}' (\ovl{0.a_1 a_2 \cdots}) = i$ if and only if $a_1 = 0$ and $t'_l(\ovl{0.a_2 a_3 \cdots}) = i - (2l + 1)h_k -l$, or $a_1 = 1$ and $t'_l(\ovl{0.a_2 a_3 \cdots}) = i - (2l + 1)h_k -l - 1$, or $a_1 = 2$ and $t'_{l+1}(\ovl{0.a_2 a_3 \cdots}) = i - 2lh_k -l$. Lastly, $t_{3l+2}' (\ovl{0.a_1 a_2 \cdots}) = i$ if and only if $a_1 = 0$ and $t'_l(\ovl{0.a_2 a_3 \cdots}) = i - (2l + 2)h_k -l - 1$, or $a_1 = 1$ and $t'_{l+1}(\ovl{0.a_2 a_3 \cdots}) = i - (2l + 1)h_k -l - 1$, or $a_1 = 2$ and $t'_{l+1}(\ovl{0.a_2 a_3 \cdots}) = i - (2l + 1)h_k -l$.
\end{proof}

Recall that $d_l(n)=\mu(t_l^{-1}(n))$. Since $d_l$ and $d_l'$ are scalar multiples of each other, the following is immediate.

\begin{corollary}
 $$ d_{3l}(i) = d_l (i - 2lh_k - l)$$
 $$ d_{3l + 1}(i)=
 \frac{1}{3} \left( d_l (i - (2l + 1)h_k - l) + d_l(i - (2l + 1)h_k - l - 1) + d_{l+1}(i - 2lh_k - l) \right)$$
 $$ d_{3l + 2}(i)=
 \frac{1}{3} \left( d_l (i - (2l + 2)h_k - l - 1) + d_{l+1}(i - (2l + 1)h_k - l - 1) + d_{l+1}(i - (2l + 1)h_k - l) \right)$$
 \label{cor:d_l}
\end{corollary}

The remainder of this paper is devoted to studying properties of $D_l$ and using them to prove our main theorems.
\begin{corollary}
 Each $D_l$ is an even function except for $x\in (1/2)\mb{Z}$, and it is increasing on $(-\infty, 0)$ and decreasing on $(0, \infty)$. Note that $\int_{-\infty}^\infty D_l(x)dx=1$. Further, $D_0 = \mathbf{1}_{[-1/2, 1/2)}$ and the following relations hold.
 
 $$ D_{3l}(x) = D_l (x)$$
 $$ D_{3l + 1}(x)=
 \frac{1}{3} \left( D_{l+1} (x) + D_l(x - 1/2) + D_l(x + 1/2) \right)$$
 $$ D_{3l + 2}(x)=
 \frac{1}{3} \left( D_l (x) + D_{l+1}(x - 1/2) + D_{l+1}(x + 1/2) \right)$$
 
 \label{cor:D_l}
\end{corollary}

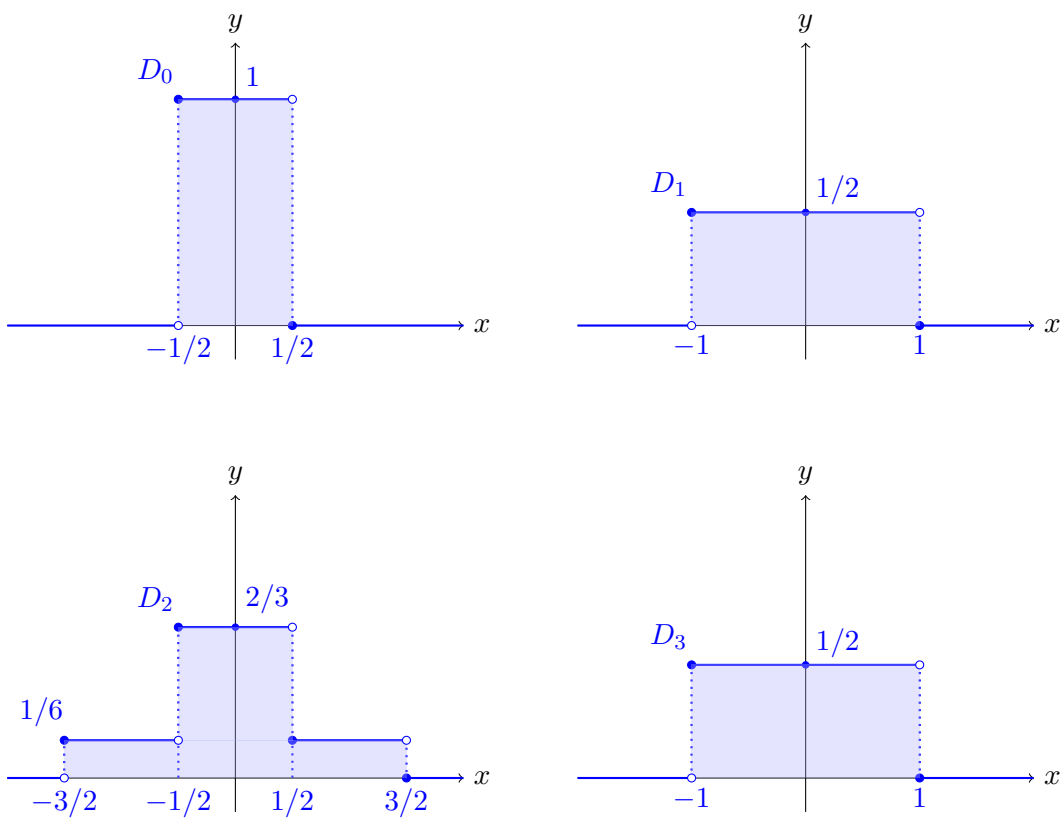
\begin{figure} [h]
\centering
 \begin{tikzpicture} [xscale = 1.5, yscale = 1.5,
 circ/.style={draw,circle,inner sep=0pt,minimum size=4pt},
 dot/.style={circ,fill}
 ]
 \draw[->] (-2,0)--(2,0) node[right]{$x$};
 \draw[->] (0,-0.3)--(0,2.5) node[above]{$y$};
 \draw[blue, thick] (-2,0)--(-0.5,0);
 
 \draw[blue, fill=blue] (-0.5,2) circle (1pt);
 \draw[blue, thick] (-0.5,2)--(0.5,2);
 
 \draw[blue, fill=blue] (0.5,0) circle (1pt)node[below]{$1/2$};;
 \draw[blue, thick] (2,0)--(0.5,0);
 \draw[blue, dotted, thick] (-0.5, 0)--(-0.5, 2);
 \draw[blue, dotted, thick] (0.5, 0)--(0.5, 2);
 \draw[blue, fill=blue] (0,2) circle (0.8pt) ;
 \node[blue,right] at (0, 2.2) {$1$};
 \path [fill=blue, fill=blue!30, opacity=0.35] (-0.5,0) rectangle (0.5,2);
 \node[blue] at (-0.7, 2.25) {$D_0$};
 \draw[blue, fill=white] (-0.5,0) circle (1pt) node[below]{$-1/2$};
 \draw[blue, fill=white] (0.5,2) circle (1pt);
 
 \draw[->] (3,0)--(7,0) node[right]{$x$};
 \draw[->] (5,-0.3)--(5,2.5) node[above]{$y$};
 \draw[blue, thick] (3,0)--(4,0);
 
 \draw[blue, fill=blue] (4,1) circle (1pt);
 \draw[blue, thick] (4,1)--(6,1);
 
 \draw[blue, fill=blue] (6,0) circle (1pt)node[below]{$1$};;
 \draw[blue, thick] (6,0)--(7,0);
 \draw[blue, dotted, thick] (4, 0)--(4, 1);
 \draw[blue, dotted, thick] (6, 0)--(6, 1);
 \draw[blue, fill=blue] (5,1) circle (0.8pt) ;
 \node[blue,right] at (5, 1.2) {$1/2$}; 
 \path [fill=blue, fill=blue!30, opacity=0.35] (4,0) rectangle (6,1);
 \node[blue] at (3.8, 1.25) {$D_1$};
 
 \draw[blue, fill=white] (4, 0) circle (1pt) node[below]{$-1$};
 \draw[blue, fill=white] (6,1) circle (1pt);
 
 \draw[->] (3,-4)--(7,-4) node[right]{$x$};
 \draw[->] (5,-4.3)--(5,-1.5) node[above]{$y$};
 \draw[blue, thick] (3,-4)--(4,-4);
 
 \draw[blue, fill=blue] (4,-3) circle (1pt);
 \draw[blue, thick] (4,-3)--(6,-3);
 
 \draw[blue, fill=blue] (6,-4) circle (1pt)node[below]{$1$};;
 \draw[blue, thick] (6,-4)--(7,-4);
 \draw[blue, dotted, thick] (4, -4)--(4, -3);
 \draw[blue, dotted, thick] (6, -4)--(6, -3);
 \draw[blue, fill=blue] (5,-3) circle (0.8pt) ;
 \node[blue,right] at (5, -2.8) {$1/2$}; 
 \path [fill=blue, fill=blue!30, opacity=0.35] (4,-4) rectangle (6,-3);
 \node[blue] at (3.8, -2.75) {$D_3$};
 
 \draw[blue, fill=white] (4, -4) circle (1pt) node[below]{$-1$};
 \draw[blue, fill=white] (6,-3) circle (1pt);
 
 \draw[->] (-2,-4)--(2,-4) node[right]{$x$};
 \draw[->] (0,-4.3)--(0,-1.5) node[above]{$y$};
 \draw[blue, thick] (-2,-4)--(-1.5,-4);
 
 \draw[blue, fill=blue] (-1.5,-3.6666) circle (1pt);
 \draw[blue, thick] (-1.5,-3.6666)--(-0.5,-3.6666);
 
 \draw[blue, fill=blue] (-0.5,-2.6666) circle (1pt);
 \draw[blue, thick] (-0.5, -2.6666)--(0.5, -2.6666);
 
 \draw[blue, thick] (0.5, -3.6666)--(1.5, -3.6666);
 \draw[blue, fill=blue] (0.5,-3.6666) circle (1pt);
 
 \draw[blue, fill=blue] (1.5,-4) circle (1pt);
 \draw[blue, thick] (2, -4)--(1.5,-4) node[below]{$3/2$};
 \draw[blue, fill=blue] (0, -2.6666) circle (0.8pt); 
 
 \node[blue,right] at (0, -2.6666+0.25) {$2/3$};
 
 \draw[blue, dotted, thick] (-1.5, -4)--(-1.5, -3.6666);
 \draw[blue, dotted, thick] (-0.5, -2.6666)--(-0.5, -4) node[below] {$-1/2$};
 \draw[blue, dotted, thick] (0.5, -2.6666)--(0.5, -4) node[below] {$1/2 $};
 \draw[blue, dotted, thick] (1.5, -3.6666)--(1.5, -4);
 
 \path [fill=blue, fill=blue!30, opacity=0.35] (-1.5,-4) rectangle (1.5,-3.6666);
 \path [fill=blue, fill=blue!30, opacity=0.35] (-0.5,-3.6666) rectangle (0.5,-2.6666);
 \node[blue] at (-0.7, -2.6666+0.25) {$D_2$};
 \node[blue] at (-1.7, -3.6666+0.25) {$1/6$};
 
 \draw[blue, fill=white] (-1.5,-4) circle (1pt) node[below]{$-3/2$}; 
 \draw[blue, fill=white] (-0.5,-3.6666) circle (1pt);
 \draw[blue, fill=white] (0.5,-2.6666) circle (1pt);
 \draw[blue, fill=white] (1.5,-3.6666) circle (1pt);
\end{tikzpicture} 
\caption{$D_l$ for $l=0,1,2,3$}
\label{fig3}
 \end{figure}

\begin{proof}
 The relations can be proven directly using Lemma~\ref{lem:d_l'}.
 
 Hence, we only need to show that each $D_l$ is an even function that increases for $x < 0$ and decreases for $x > 0$. Basic calculations show the claim for $l=0,1,2,3$. The induction uses the elementary fact that finite sums and the operator $f\mapsto f(\cdot-1/2)+f(\cdot+1/2)$ preserve even unimodality for step functions on the $\frac12\mathbb Z$-grid.
\end{proof}
Let us calculate $d_l'$ for $l=0,1,2,3$.
\begin{enumerate}
 \item $l=0$ : $d_0'(0)=1$
 \item $l=1$ : $d_1'(h_k)=1/2$, $d_1'(h_k+1)=1/2$
 \item $l=2$ : $d_2'(2h_k)=1/6$, $d_2'(2h_k+1)=2/3$, $d_2'(2h_k+2)=1/6$
 \item $l=3$ : $d_3'(3h_k+1)=1/2$, $d_3'(3h_k+2)=1/2$
\end{enumerate}
Note that the functions $d_l'$ are zero elsewhere. We have $D_l$ for $l=0,1,2,3$ as illustrated in Figure~\ref{fig3}.

\subsection{\texorpdfstring{Properties of \(D_l\)}{Properties of Dl}}
\label{subsec:D_l-properties}

Now we use Corollary~\ref{cor:D_l} to prove properties of $D_l$. To this end, define
\[ 
B_l := \{ n \in \mathbb{Z} : d_l (n) > 0 \}
\]
and $b_l = |B_l|$. Corollary~\ref{cor:D_l} tells us that $B_l$ is a contiguous subset of the integers. The following lemmas give a recursive equation for $b_l$ and show that the height of $D_l$ and $\| D_l - D_{l+1} \|_1$ are both bounded by $O(b_l^{-1/2})$. Since $D_l$ is unimodal, this also implies that $V(D_l)$ decays as $b_l$ diverges.

\begin{lemma}
 \label{lem:b_l-recursion}
 $b_l$ is given by $b_0 = 1$, $b_1=2$, and
 \[ 
 b_{3l} = b_l, \quad b_{3l+1} = \max(b_l + 1 ,\, b_{l+1}), \quad b_{3l+2} = \max(b_l,\, b_{l+1} + 1).
 \]
 Furthermore, this implies the following properties of $b_l$.
 \begin{enumerate}
 \item $|b_l - b_{l+1}| = 1$.
 \item $b_{3l+1} = b_l + 1$ for $l\ge0$, and $b_{3l-1}=b_l+1$ for $l\ge1$.
 \end{enumerate}
\end{lemma}
\begin{proof}
Define $B'_l=\{x\in \mb{R}\,:\, D_l(x)> 0\}$. By Corollary~\ref{cor:D_l}, we see that $B'_l$ is an interval centered at the origin, $b_l = |B'_l|$, and
$$B'_{3l} = B'_l$$
$$B'_{3l+1} = B'_{l+1} \cup \left(B'_l - \frac{1}{2}\right) \cup \left(B'_l + \frac{1}{2}\right)$$
$$B'_{3l+2} = B'_{l} \cup \left(B'_{l+1} - \frac{1}{2}\right) \cup \left(B'_{l+1} + \frac{1}{2}\right).$$
These facts imply our claim. Note that the claim gives a recursive description of $(b_l)$ once $b_0 = 1$ and $b_1=2$ are fixed. The properties are easy to check.
\end{proof}

 \begin{lemma}
 Let $H_l=\max_{x\in \mb{R}} D_l(x)$. Then,
 \begin{equation*}
 H_l \le \frac{C}{b_l^{1/2}}.
 \label{eq:H_l}
 \end{equation*}
 \label{lem:H_l}
 \end{lemma}
 \begin{proof}
 Note that $H_l = D_l(0)$. We prove this lemma by comparing $D_l$ to the distribution of a lazy random walk on $\mathbb{Z}$. The local limit theorem gives a bound on the lazy random walk, which completes the proof.
 \begin{proof}[Step 1]\renewcommand{\qedsymbol}{}
 Define $\phi : L^1(\mathbb{R}) \to L^1(\mathbb{R})$ so that $$\phi(f)(x) = \frac{1}{2}(f(x - 1/2) + f(x + 1/2)).$$Corollary~\ref{cor:D_l} gives $$D_{3l\pm1}=\dfrac{1}{3}\bkt{D_l(\cdot - 1/2)+D_l(\cdot + 1/2)}+\dfrac{1}{3}D_{l\pm1}=\dfrac{2}{3}\phi D_l+\dfrac{1}{3}D_{l\pm 1}.$$ By induction, we check that every $D_l$ can be represented as $$D_l=(c_{b_l-1}\phi^{b_l-1}+\cdot +c_1\phi+c_0)D_0$$ for some $\sum_{i=0}^{b_l-1} c_i=1$ with $c_{b_l-1}\neq 0$. It is also routine to check that such a representation is unique.
 \end{proof}
 
 \begin{proof}[Step 2]\renewcommand{\qedsymbol}{}
 For two polynomials $f(x)=c_mx^m+\cdots +c_1x+c_0$ and $g(x)=c_m'x^m+\cdots +c_1'x+c_0'$, we say $f \preceq g$ if $\sum_{i=0}^t c_i\leq \sum_{i=0}^t c_i'$ for every $0\leq t\leq m$ (it is clear that $\preceq$ is a partial ordering). Note that $c_m$ or $c_m'$ need not be nonzero. The following are routine to check.
 \begin{enumerate}
 \item if $f\preceq g$, then $pf\preceq pg$ for any polynomial $p$ with positive coefficients.
 \item if $f_1 \preceq g_1$ and $f_2\preceq g_2$, then $f_1+f_2\preceq g_1+g_2$.
 \item if $f \preceq g$, then $f(\phi)(D_0)(0) \le g(\phi)(D_0)(0)$.
 \end{enumerate}
 We also use this notation to compare functions of the form $f(\phi)(D_0)$. For instance, since $D_1 = \phi(D_0)$, we can write $D_1 \preceq D_0$.
 \end{proof}
 
 \begin{proof}[Step 3]\renewcommand{\qedsymbol}{}
 Let $F_0=D_0$. We define $$F_n=\bkt{\dfrac{1}{3}\phi+\dfrac{2}{3}}F_{n-1}$$ for $n\geq 1$. We claim the following three facts:
 \begin{enumerate}
 \item $F_n \preceq F_{n-1}$.
 \item $\frac{2}{3}\phi F_{n-1} + \frac{1}{3}F_{n-2} \preceq F_n$.
 \item $D_l \preceq F_{b_l-1}$.
 \end{enumerate}
 (1) is clear from the definition of $F_n$, along with the fact that $\frac{1}{3} \phi + \frac{2}{3} \preceq 1$. For (2), we have $$\dfrac{2}{3}\phi F_{n-1}+\dfrac{1}{3}F_{n-2}=\bkt{\dfrac{2}{9}\phi^2+\dfrac{4}{9}\phi+\dfrac{1}{3}}F_{n-2}\preceq\bkt{\dfrac{1}{9}\phi^2+\dfrac{4}{9}\phi+\dfrac{4}{9}}F_{n-2}=F_n.$$
 To show (3), we use induction on $l$. When $l = 0, 1$, the claim is easy to check. For larger values, note that
 $$D_{3l} = D_{l} \preceq F_{b_l - 1} = F_{b_{3l} - 1}$$
 and 
 $$D_{3l \pm 1} = \frac{2}{3}\phi D_l + \frac{1}{3}D_{l \pm 1} \preceq \frac{2}{3}\phi F_{b_l - 1} + \frac{1}{3}F_{b_{l \pm 1} - 1} \preceq \frac{2}{3}\phi F_{b_l - 1} + \frac{1}{3} F_{b_l - 2} \preceq F_{b_l} = F_{b_{3l \pm 1} - 1}$$
 (note that $b_l=b_{3l}=b_{3l\pm1}-1$ by Lemma~\ref{lem:b_l-recursion}, and $b_{l \pm 1} - 1 \ge b_l - 2$). As such, we can use induction to show that (3) holds for all $l$.
 \end{proof}
 
 \begin{proof}[Step 4]\renewcommand{\qedsymbol}{}
 Note that $$F_n(x) = \frac{1}{3}\phi F_{n-1}(x) + \frac{2}{3}F_{n-1}(x) = \frac{1}{6}F_{n-1}(x-1/2) + \frac{1}{6}F_{n-1}(x+1/2) + \frac{2}{3}F_{n-1}(x).$$
 As such, $F_n(0)$ is bounded by the sum of at most two adjacent point probabilities for the lazy random walk on $\frac12\mathbb Z$ with $n$ steps that moves to the left or right by $1/2$ with probability $1/6$ and stays put with probability $2/3$. By the local limit theorem \cite{Durrett}, we see that $F_{n}(0) = O(n^{-1/2})$. Together with the fact that $D_l \preceq F_{b_l - 1}$, our proof is complete.
 \end{proof}
\end{proof}

\begin{lemma}
\label{lem:Dl-Dlp1}
 \[
 \| D_l - D_{l+1} \|_1 \le \frac{C}{b_l^{1/2}}.
 \]
\end{lemma}

\begin{proof}
 We use induction on $l$. The base case is trivial. Assume as the induction hypothesis that $\| D_l - D_{l + 1 } \|_1 \le Cb_{l}^{-1/2}$. Since $D_l$ is unimodal, $\|D_l-D_l(\cdot \pm 1/2)\|_1\le H_l$. Then
 \begin{align*}
 \|D_{3l} - D_{3l + 1}\|_1 &\le \frac{1}{3} \left\{ \|D_l - D_{l + 1}\|_1 + \|D_l - D_l(\cdot - 1/2)\|_1 + \|D_l - D_l(\cdot + 1/2)\|_1 \right\} \\
 &\le \frac{1}{3} \left( \frac{C+2C'}{b_l^{1/2}} \right)
 \end{align*}
 by Lemma~\ref{lem:H_l} and the induction hypothesis. Here, $C'$ is the constant given by Lemma~\ref{lem:H_l}. Since $b_l = b_{3l}$, we can choose a sufficiently large $C$ to complete the induction. The other cases follow similarly, sometimes utilizing the fact that $|b_l - b_{l+1}| = 1$.
\end{proof}

We conclude with interesting characterization for $b_l$, relating it to the \emph{balanced ternary expansion} of $l$. This result is not strictly necessary since we prove a more general result in the next section. Still, we include it here as an interesting fact.

\begin{lemma}
 For every $l\in \mb{N}_{\ge 0}$, there exists a unique sequence $a_{l,i}\in \{-1, 0, 1\}$, with finitely many nonzero terms, such that 
 \begin{equation}
 l=\sum_{i=0}^\infty a_{l,i} \,3^i.
 \label{eq:ternary}
 \end{equation}
 This expansion is called the \emph{balanced ternary expansion} of $l$. Then we obtain $$b_l=1+\sum_{i=0}^\infty \abs{a_{l,i}}.$$
 \label{lem:ternary}
\end{lemma}
\begin{proof}
 It is routine to check that $a_{l, i}$ exists uniquely (\eqref{eq:ternary} is called the \emph{balanced ternary system}). Let $$b_l' = 1 + \sum_{i=0}^{\infty} | a_{l, i}|.$$ We show that $b_l'$ also satisfies the recursive formula of Lemma~\ref{lem:b_l-recursion}. Observe that
 $$ a_{3l, i} = 
 \begin{cases}
 a_{l, i-1} & \textrm{if } i > 0 \\
 0 & \textrm{if } i = 0\\
 \end{cases}$$
 $$ a_{3l + 1, i} = 
 \begin{cases}
 a_{l, i-1} & \textrm{if } i > 0 \\
 1 & \textrm{if } i = 0\\
 \end{cases}$$
 $$ a_{3l + 2, i} = 
 \begin{cases}
 a_{l+1, i-1} & \textrm{if } i > 0 \\
 -1 & \textrm{if } i = 0.\\
 \end{cases}$$
 for any nonnegative integer $l$. Thus, the following hold:
 \begin{enumerate}
 \item $b'_{3l} = b'_l$
 \item $b'_{3l+1} = b'_l + 1$
 \item $b'_{3l+2} = b'_{l+1} + 1$
 \item $|b'_{3l} - b'_{3l+1}| = 1$
 \item $|b'_{3l+1} - b'_{3l+2}| = |b'_l - b'_{l+1}|$
 \item $|b'_{3l+2} - b'_{3l+3}| = 1.$
 \end{enumerate}
 Since $|b'_0 - b'_1| = |1 - 2| = 1$, the last three properties show (recursively) that $|b'_l - b'_{l+1}| = 1$ for all $l$. This shows that the first three properties imply the recursive formula stated above (since $b'_{l+1} \le b'_l + 1$ and $b'_l \le b'_{l+1} + 1$), so the proof is complete.
\end{proof}

\begin{corollary}
 For any fixed $t > 0$, we obtain
 \[
 | \{ l \le n : b_l \le t \}| \le (C \log n)^t.
 \]
\end{corollary}
\begin{proof}
Since \(l \le n\), the balanced ternary expansion of \(l\) (i.e. the digits \(a_{l,i}\)) can only contain nonzero digits within the first \(O(\log n)\) positions. Considering all numbers of at most $C\log n$ nonzero digits in the balanced ternary system, the desired results follows from Lemma~\ref{lem:ternary} and standard binomial coefficient estimates.
\end{proof}

\section{Tight Maps}
\label{sec:rank-one}

\subsection{Rank-one transformations}
We now introduce rank-one transformations. There are several equivalent definitions; here we present Definition~4 of~\cite{Ferenczi}. Rank-one transformations generalize the Chacon transformation by allowing the number of towers and spacers to change at each step. They are uniquely determined by the doubly-indexed \emph{spacer sequence} $(s_{n, j} , m_n)$. In this sequence, $m_n$ represents the number of towers at step $n$, and $s_{n, j}$ denotes the number of spacers placed above the $j$-th tower in that step, where $0 \le j < m_n$. The \emph{cut sequence} $\{ m_n \}$ determines the number of towers at each step, and the \emph{height sequence} $\{ h_n \}$ is given by $h_0 = 1$ and $h_{n+1} = m_n h_n + \sum_{j = 0}^{m_n - 1} s_{n, j}$. Note that the spacer sequence completely determines the lengths of each tower and spacer. For more information and properties of rank-one transformations, see \cite{Ferenczi}, \cite{Friedman2}, and \cite{Adams}. For example, the Chacon transformation is a rank-one transformation with $m_n = 3$ and spacer sequence $(0,1,0)$.

In this paper, we will study rank-one transformations with constant spacer sequences, which we call \emph{tight maps}. Specifically, we will construct universal exceptional sets for an even smaller subclass of \emph{restrictive} tight maps. The lower-bound argument will only use the assumption that no spacers are placed above the last subcolumn, namely $s_{m-1}=0$. Clearly, the Chacon transformation is an example of a restrictive tight map.

Throughout, rank-one systems are understood in the standard finite-measure cutting-and-stacking model: the tower levels generate the ambient sigma-algebra modulo null sets, and the tower unions exhaust the space modulo null sets. Equivalently,
\(\mcr B=\sigma\{T^rA_k:k\ge0,\ 0\le r<h_k\}\) modulo null sets.

\begin{definition}
 A rank-one transformation is called \textit{tight} if its cutting sequence is the same at each stage and its spacer sequence is constant regardless of \(n\). That is, \(m_n = m\) is fixed and the spacer sequence can be written as \(s_{n,j} = s_j\).
 \label{def1}
\end{definition}

For tight maps, we use the convention \(s_{j+m}=s_j\) for \(j\in\mathbb Z\).

\begin{definition}
 A tight map with spacer sequence $(s_0 , s_1 , \dots , s_{m-1})$ is \textit{restrictive} if it has the following properties.
 \begin{enumerate}
 \item $s_{m-1} = 0$.
 \item $\gcd(s_0 , \dots , s_{m-1}) = 1$.
 \item For every $1\le r<m$, each of the following two integers is either $0$ or $1$, and at least one of them is $1$:
 \[
 \gcd\{s_{j+r}-s_j:0\le j\le m-r-2\},
 \]
 \[
 \gcd\{s_{j+r-m}-s_j:m-r\le j\le m-2\}.
 \]
 Here we use the convention that the gcd of an empty set, as well as the gcd of \(\{0\}\), is \(0\).
 \end{enumerate}
 \label{def2}
\end{definition}

For tight maps, we can repeat Section~\ref{subsec:chacon-recursion} to get a recursive formula for $D_l$. Let $A_k$ be the bottom interval in the $k$-th step of the cutting and stacking process, and define $a_k = \mu(A_k)^{-1}$, $r_k$, $S_k$, $t_l'$, and $d_l'$ as before. Let $s := \sum_{j=0}^{m-1} s_j$, and numbers with an overline denote numbers expressed in base-$m$. Since $h_{k+1}=mh_k+s$ and $h_0=1$, we have $h_k=m^k+s(m^k-1)/(m-1)$; after normalizing the total measure to one,
\[
 a_k=m^k\left(1+\frac{s}{m-1}\right)=h_k+\frac{s}{m-1}.
\]

\begin{lemma}
Given a tight map $T$, we have
 $$ r_k(\ovl{0.a_1 a_2 a_3 \cdots})=
 \begin{cases}
 h_k + s_{a_1} & \textrm{if } a_1 \ne m-1 \\
 r_k(\ovl{0.a_2 a_3 \cdots}) + s_{a_1}& \textrm{if } a_1 = m-1,
 \end{cases}$$
 $$ S(\ovl{0.a_1 a_2 a_3 \cdots})=
 \begin{cases}
 \ovl{0.(a_1+1) a_2 a_3 \cdots} & \textrm{if } a_1 \ne m-1 \\
 \frac{1}{m}S(\ovl{0.a_2 a_3 \cdots}) & \textrm{if } a_1 = m-1.
 \end{cases}$$
 \label{lem:r_k-gen}
\end{lemma}

\begin{proof}
 The proof is similar to Lemma~\ref{lem:r_k}. For any $a_1$, $T^{h_k + s_{a_1} - 1}$ is at the top of the tower. If $a_1 < m-1$, this implies that at the next turn, it will get mapped back to $A_k$, shifting by one column as it does so.

 If $a_1 = m-1$, then $T^{h_k + s_{a_1} - 1}$ is at the top of the tower even after the stacking operation. Its relative position in the tower is $\ovl{0.a_2 a_3 \dots}$, so it takes $r_k(\ovl{0.a_2 a_3 \dots}) - h_k + 1$ additional turns to get mapped back to $A_k$, in which case it gets mapped to $\frac{1}{m}S(\ovl{0.a_2 a_3 \cdots})$.
\end{proof}

\begin{remark}
$S$ depends only on $m$. That is, it is independent of both $k$ and the spacer sequence $\{s_j\}$.
\end{remark}

\begin{corollary}
\label{cor:S-gen}
Suppose \(T\) is a tight map. Then,
\[
S^{l}\!\left(\overline{0.a_1a_2a_3\cdots}\right)
=
\frac1m
\left(
(a_1+l)-m\left\lfloor\frac{a_1+l}{m}\right\rfloor
+
S^{\left\lfloor\frac{a_1+l}{m}\right\rfloor}
\!\left(\overline{0.a_2a_3\cdots}\right)
\right).
\]
Equivalently,
\[
S^{l}\!\left(\overline{0.a_1a_2a_3\cdots}\right)
=
\frac1m
\left(
(a_1+l)\bmod m
+
S^{\left\lfloor\frac{a_1+l}{m}\right\rfloor}
\!\left(\overline{0.a_2a_3\cdots}\right)
\right).
\]
\end{corollary}

\begin{proof}
For \(0\le l<m\), this follows directly from Lemma~\ref{lem:r_k-gen}. The general case follows by iterating the same carry rule in base \(m\).
\end{proof}

\begin{lemma}
 For a tight map $T$, we have $t_0' = 0$ and
 \[ 
 t_{ml + r}'(\ovl{0.a_1 a_2 a_3 \dots}) = 
 \begin{cases}
 (ml - l + r)h_k + sl + \sum_{i=0}^{r-1}s_{a_1 + i} + t_l'(\ovl{0.a_2 a_3 \dots}) & \text{ if } a_1 + r < m \\
 (ml - l + r - 1)h_k + sl + \sum_{i=0}^{r-1}s_{a_1 + i} + t_{l+1}'(\ovl{0.a_2 a_3 \dots}) & \text{ if } a_1 + r \ge m.
 \end{cases}
 \]
 \label{lem:t_l'-gen}
\end{lemma}

\begin{proof}
First note that
\begin{align*}
r_k(S^i(\ovl{0.a_1 a_2 \dots})) &= r_k \left( \frac{1}{m} \left\{ (a_1+i)-m\left\lfloor \frac{a_1+i}{m}\right\rfloor + S^{\left\lfloor \frac{a_1 + i}{m} \right\rfloor} ( \ovl{0. a_2 a_3 \cdots}) \right\} \right) \\
&= \begin{cases}
 h_k + s_{a_1 + i} & \text{ if } a_1 + i \not \equiv m-1 \mod m \\
 s_{a_1 + i} + r\left(S^{\left\lfloor \frac{a_1 + i}{m} \right\rfloor} ( \ovl{0. a_2 a_3 \cdots}) \right) & \text{ if } a_1 + i \equiv m-1 \mod m.
\end{cases}
\end{align*}
By this convention, we see that

 \begin{align*}
 t_{ml + r}'(\ovl{0.a_1 a_2 a_3 \dots}) &= \sum_{i = 0}^{ml + r - 1} r_k (S^i(\ovl{0.a_1 a_2 a_3 \dots})) \\
 &= \begin{cases}
 (ml - l + r)h_k + \sum_{i=0}^{ml+r-1}s_{a_1 + i} + t_l'(\ovl{0.a_2 a_3 \dots}) & a_1 + r < m \\
 (ml - l + r - 1)h_k + \sum_{i=0}^{ml+r-1}s_{a_1 + i} + t_{l+1}'(\ovl{0.a_2 a_3 \dots}) & a_1 + r \ge m
 \end{cases} \\
 &= \begin{cases}
 (ml - l + r)h_k + sl + \sum_{i=0}^{r-1}s_{a_1 + i} + t_l'(\ovl{0.a_2 a_3 \dots}) & a_1 + r < m \\
 (ml - l + r - 1)h_k + sl + \sum_{i=0}^{r-1}s_{a_1 + i} + t_{l+1}'(\ovl{0.a_2 a_3 \dots}) & a_1 + r \ge m.
 \end{cases} 
 \end{align*}
\end{proof}

Note that this lemma implies the following equation. Let
\[
I_l:=\int_0^1 t_l'(x)\,dx.
\]
Averaging Lemma~\ref{lem:t_l'-gen} over the first digit gives
\[
I_{ml+r}
=
\left(ml-l+r-\frac rm\right)h_k
+sl+
\frac{rs}{m}
+\frac{m-r}{m}I_l+
\frac rm I_{l+1}.
\]
By induction, this implies that
\[
I_l=\left(h_k+\frac{s}{m-1}\right)l=a_kl=\mu(A_k)^{-1}l.
\]
Therefore,
\begin{equation}
 \int_{A_k} t_l\,d\mu = l
 \label{eq:kac-lemma}
\end{equation}
even when \(T\) is not ergodic.

\begin{lemma}
 Given a tight map $T$, we have $d_0' = \mathbf{1}_{\{ 0 \} }$ and
 \begin{align*}
 d_{ml + r}'(i) = & \frac{1}{m} \sum_{j=0}^{m-r-1} d_l'\left(i - (ml - l + r)h_k -sl - \sum_{u = 0}^{r-1} s_{j+u}\right) \\
 &+ \frac{1}{m} \sum_{j=m-r}^{m-1} d_{l+1}'\left(i - (ml - l + r - 1)h_k -sl - \sum_{u = 0}^{r-1} s_{j+u}\right).
 \end{align*} 
 \label{lem:d_l'-gen}
\end{lemma}

\begin{proof}
This proof is similar to that of Lemma~\ref{lem:d_l'}. 
\end{proof}

Recall that $d_l(n)=\mu(t_l^{-1}(n))$. Since $d_l$ and $d_l'$ are scalar multiples of each other, the following is immediate.

\begin{corollary}
\label{cor:d_l-gen}
When $T$ is a tight map, we have $d_0 = \mu(A_k) \mathbf{1}_{\{ 0 \}}$ and
 \begin{align*}
 d_{ml + r}(i) = & \frac{1}{m} \sum_{j=0}^{m-r-1} d_l\left(i - (ml - l + r)h_k -sl - \sum_{u = 0}^{r-1} s_{j+u}\right) \\
 &+ \frac{1}{m} \sum_{j=m-r}^{m-1} d_{l+1}\left(i - (ml - l + r - 1)h_k -sl - \sum_{u = 0}^{r-1} s_{j+u}\right).
 \end{align*} 
\end{corollary}

From now on, we will use convolution by distributions $\alpha_r$ and $\beta_r$ to denote these relations.
\begin{corollary}
Suppose $T$ is a tight map. Then,
 \[
 D_{ml + r}(x) = \frac{1}{m} \left[ \sum_{j=0}^{m-r-1} D_l \left( x + \frac{sr}{m-1} - \sum_{u=0}^{r-1} s_{j+u} \right) + \sum_{j=m-r}^{m-1} D_{l+1} \left( x + \frac{s(r-1)}{m-1} - \sum_{u=0}^{r-1} s_{j+u} \right) \right].
 \]
 \label{cor:D_l-gen}
 In other words,
 \begin{equation}
 \label{eq:recursion}
 D_{ml + r} = \frac{m-r}{m} \alpha_r \ast D_l + \frac{r}{m} \beta_r \ast D_{l+1},
 \end{equation}
 where for $1\le r<m$, $\alpha_r$, $\beta_r$ are probability distributions
 \[
 \alpha_r = \frac{1}{m-r} \left( \sum_{j=0}^{m-r-1} \delta \left( \sum_{u=0}^{r-1} s_{j+u} - \frac{sr}{m-1} \right) \right), \quad \beta_r = \frac{1}{r} \left( \sum_{j=m-r}^{m-1} \delta \left( \sum_{u=0}^{r-1} s_{j+u} - \frac{s(r-1)}{m-1} \right) \right),
 \]
with the conventions $\alpha_0=\beta_0=\alpha_m=\beta_m=\delta_0$ whenever these symbols are used.
\end{corollary}

\begin{proof}
 The relations can be proven directly using Corollary~\ref{cor:d_l-gen}.
\end{proof}

\begin{remark}
Condition (3) in Definition~\ref{def2} says precisely that, for every
$1\le r<m$, the support differences of each of $\alpha_r$ and $\beta_r$
generate either $\{0\}$ or $\mathbb Z$, and at least one of them generates
$\mathbb Z$.
\end{remark}

\subsection{\texorpdfstring{Properties of \(D_l\)}{Properties of Dl}}

In this section, we bound the support of $D_l$ and measure the decay of $V(D_l)$. The support estimate only uses the assumption that no spacers are placed above the last subcolumn, while the variation estimates use the full restrictive hypotheses.

\begin{lemma}
The following are true for restrictive tight maps.
 \begin{enumerate}
 \item $D_{ml}=D_l$.
 \item For any $r \ne 0$, the support differences of each of $\alpha_r$ and $\beta_r$ generate either $\{0\}$ or $\mathbb Z$, and at least one of them generates $\mathbb Z$. In particular, $\alpha_r$ and $\beta_r$ are not both Dirac distributions.
 \end{enumerate}
 \label{lem:reduced-prop}
\end{lemma}

\begin{proof}
\quad
 \begin{enumerate}
 \item This is the case $r=0$ of Corollary~\ref{cor:D_l-gen}, since $\alpha_0=\beta_0=\delta_0$ by convention.
 \item This is exactly condition (3) in Definition~\ref{def2}, translated through the definitions of $\alpha_r$ and $\beta_r$ in Corollary~\ref{cor:D_l-gen}. The final assertion follows because a Dirac distribution has support-difference group $\{0\}$.
 \end{enumerate}
\end{proof}

We define the following sequences $b_n$ and $c_n$. Note that $b_n$ bounds the support of $D_l$ (Lemma~\ref{lem:supp-bound}), while $c_l$ measures the amount of convolutions applied to $D_l$. Thus, $c_l$ is related to $V(D_l)$ under the restrictive hypotheses (Lemma~\ref{lem:variation-bound}).
\begin{align}
\label{eq:bn-def}
b_0 = 1, \quad b_1 = 2 , \quad b_{ml + r} &= \begin{cases}
 b_l & r = 0 \\
 \max(b_l + 1, b_{l+1}) & r = 1 \\
 \max(b_l, b_{l+1}) + 1 & 2 \le r \le m-2 \\
 \max(b_l, b_{l+1} + 1) & r = m-1
\end{cases} \\
\label{eq:cn-def}
c_0 = 1, \quad c_1 = 2 , \quad c_{ml + r} &= \begin{cases}
 c_l & r = 0 \\
 \min(c_l, c_{l+1}) + 1 & 1 \le r < m
\end{cases}
\end{align}

\begin{lemma}\label{lem:c-diff-bound}
For the sequence \((c_l)_{l\ge0}\) defined in \eqref{eq:cn-def},
\[
|c_l-c_{l+1}|\le1
\]
for every \(l\ge0\).
\end{lemma}
\begin{proof}
Induct on \(n=ml+r\). If \(r=0\), then \(|c_{ml+1}-c_{ml}|=|\min(c_l,c_{l+1})+1-c_l|\le1\). If \(1\le r<m-1\), the two values are equal. If \(r=m-1\), then \(|c_{ml+m}-c_{ml+m-1}|=|c_{l+1}-\min(c_l,c_{l+1})-1|\le1\).
\end{proof}

The following lemma relates the support of $D_l$ to $b_l$. We need one more auxiliary constant:
\begin{equation}
 \label{eq:R-def}
 R := \max\left\{\frac12,\ \max_r \max \{ |x| : x\in \operatorname{supp}\alpha_r \cup \operatorname{supp}\beta_r \}\right\}.
\end{equation}
\begin{lemma}
\label{lem:supp-bound}
Let $T$ be a tight map with $m\ge3$ and spacer sequence $(s_0,\ldots,s_{m-1})$, and assume that $s_{m-1}=0$. Then,
\[
\supp(D_l)\subseteq[-Rb_l,Rb_l].
\]
\end{lemma}
\begin{proof}
We argue by induction on $l$. The case $l=0$ follows from $\operatorname{supp}D_0\subseteq[-1/2,1/2]$ and the definition of $R$. For $l=1$, the case $ml+r=1$ of \eqref{eq:recursion} gives
\[
D_1=\frac{m-1}{m}\alpha_1*D_0+\frac1m\beta_1*D_1.
\]
Since $s_{m-1}=0$ gives $\beta_1=\delta_0$, solving this identity yields $D_1=\alpha_1*D_0$; hence the claim follows from the definition of $R$ and $b_1=2$.

Assume the claim for all indices smaller than $ml+r$, where $0\le r<m$ and $ml+r>1$. If $\gamma$ is one of the distributions $\alpha_r,\beta_r$, then the induction hypothesis gives
\[
\supp(\gamma*D_i)\subseteq
\begin{cases}
[-Rb_i,Rb_i],& \gamma=\delta_0,\\
[-R(b_i+1),R(b_i+1)],& \gamma\ne\delta_0.
\end{cases}
\]
Applying this to the two terms in \eqref{eq:recursion}, and using the recursive definition of $b_l$, gives the desired bound. The endpoint cases use
\[
\alpha_0=\beta_0=\delta_0
\]
by convention and
\[
\alpha_{m-1}=\beta_1=\delta_0,
\]
which follow from $s_{m-1}=0$.
\end{proof}

\begin{lemma}
\label{lem:variation-bound}
 Suppose that $\alpha$ is a probability distribution supported on a finite subset of a coset of $\mathbb Z$, and that its support differences generate $\mathbb Z$. Then, for $n\ge1$,
 \[ 
 V(\alpha^{\ast n} \ast D_0) \le \frac{C_\alpha}{\sqrt{n}},
 \]
 where $C_\alpha$ is some constant that may depend on $\alpha$.
\end{lemma}

\begin{proof}
Choose $a\in\mathbb R$ such that $\operatorname{supp}\alpha\subset a+\mathbb Z$,
and let $\widetilde\alpha$ be the translate of $\alpha$ by $-a$. Then
$\alpha^{\ast n}\ast D_0$ is a translate of
$\widetilde\alpha^{\ast n}\ast D_0$, so their total variations are equal.
Moreover, the support-difference assumption is unchanged by this translation.
Thus we may assume that $\alpha$ is supported on $\mathbb Z$.

Let $X_1 , X_2 , \dots$ be i.i.d. variables with distribution $\alpha$. Then,
$\alpha^{\ast n}$ is the probability distribution of
$S_n^X = \sum_{i=1}^n X_i$. Similarly, let $Y_1, Y_2 , \dots $ be a not
necessarily independent copy of $X_1 , X_2 , \dots$ and let
$S_n^Y = \sum_{i=1}^n Y_i - 1$. Thus, $S_n^Y$ has distribution
$\alpha^{\ast n} \ast \delta(-1)$. Since $D_0$ is uniformly distributed on
$[-1/2, 1/2)$,
 \begin{align*}
 V(\alpha^{\ast n} \ast D_0) &= \sum_{x \in \mathbb{Z}} \left|( \alpha^{\ast n} \ast D_0 )(x) - (\alpha^{\ast n} \ast D_0 )(x-1) \right| \\
 &= \sum_{x \in \mathbb{Z}} \left| \mathbb{P}(S_n^X = x) - \mathbb{P}(S_n^Y = x) \right| \\
 &\le 2 \mathbb{P}(S_n^X \ne S_n^Y).
 \end{align*}
 We couple the two walks so that, once they meet, they move together. Namely,
 if $S_{n-1}^X = S_{n-1}^Y$, choose $X_n=Y_n$ according to $\alpha$;
 otherwise choose $X_n$ and $Y_n$ independently according to $\alpha$. Let
 $\tau=\min\{n:S_n^X=S_n^Y\}$. Then
 \[
 V(\alpha^{\ast n}\ast D_0)\le 2\mathbb P(\tau>n).
 \]
 Until time $\tau$, the difference $S_n^X-S_n^Y$ is a symmetric random walk
 on $\mathbb Z$ starting from $1$, with increment distribution
 $\alpha*(-\alpha)$. By the assumption on the support of $\alpha$, this
 difference walk is irreducible; and since $0\in\operatorname{supp}(\alpha*(-\alpha))$,
 it is aperiodic. The classical one-dimensional hitting-time estimate gives
 \[
 \mathbb P(\tau>n)=O(n^{-1/2}),
 \]
 which completes the proof.
\end{proof}

\begin{lemma}
 Suppose $T$ is a restrictive tight map. Then,
 \begin{equation*}
 V(D_l) \le \dfrac{C}{c_l^{1/2}}.
 \label{eq:V(D_l)}
 \end{equation*}
 \label{lem:tot-var-bound}
\end{lemma}
\begin{proof}
Let \(\Gamma\) be the finite set of all non-Dirac distributions among
\(\alpha_r,\beta_r\), \(1\le r<m\). By Lemma~\ref{lem:reduced-prop}, every
\(\gamma\in\Gamma\) has support differences generating \(\mathbb Z\). Hence
Lemma~\ref{lem:variation-bound} gives
\[
V(\gamma^{*N}*D_0)\le C_\gamma N^{-1/2}\qquad(N\ge1),
\]
and we choose one constant valid for all \(\gamma\in\Gamma\).

Iterating \eqref{eq:recursion} expands \(D_l\) as a convex combination
\[
D_l=\sum_\omega c_\omega\;\gamma_{\omega,1}*\cdots*\gamma_{\omega,N(\omega)}*D_0,
\qquad c_\omega\ge0,\quad \sum_\omega c_\omega=1,
\]
where the \(\gamma_{\omega,i}\) are elements of \(\Gamma\), and all Dirac factors have been omitted. If \(N(\omega)\ge1\), some member of the finite set \(\Gamma\) occurs at least \(N(\omega)/|\Gamma|\) times. Since convolution is commutative and convolution by a probability distribution does not increase total variation,
\begin{equation}
\label{eq:variation-many-factors}
V(\gamma_{\omega,1}*\cdots*\gamma_{\omega,N(\omega)}*D_0)
\le C N(\omega)^{-1/2}.
\end{equation}
Thus it remains to control the total coefficient of terms with few non-Dirac factors.

Let \(F_l(z)\) be the probability polynomial whose coefficient of \(z^N\) is the total coefficient of terms with exactly \(N\) non-Dirac factors in the above expansion of \(D_l\). Put
\[
\eta_r=\mathbf1_{ \{ \alpha_r\text{ is non-Dirac}\} },
\qquad
\theta_r=\mathbf1_{ \{ \beta_r\text{ is non-Dirac}\} }.
\]
Then
\[
F_{ml}=F_l,\qquad
F_{ml+r}=\frac{m-r}{m}z^{\eta_r}F_l+\frac r m z^{\theta_r}F_{l+1}\quad(1\le r<m),
\]
with \(F_0=1\) and \(F_1=z\); the latter uses the identity \(D_1=\alpha_1*D_0\) established in the proof of Lemma~\ref{lem:supp-bound}.

For probability polynomials \(P,Q\), write \(P\preceq Q\) if every lower partial sum of the coefficients of \(P\) is bounded by the corresponding lower partial sum of \(Q\). This order is preserved by convex combinations and by multiplication by a fixed probability polynomial. Let
\[
q(z)=\frac{m-1}{m}+\frac1m z.
\]
We claim that
\begin{equation}
\label{eq:Fl-domination}
F_l\preceq q^{c_l-1}
\end{equation}
for every \(l\). The cases \(l=0,1\) are immediate. The step \(ml\) is inherited from \(l\). For \(ml+r\), \(1\le r<m\), set \(C_*=\min(c_l,c_{l+1})\). By Lemma~\ref{lem:c-diff-bound}, the exponents \(c_l-1\) and \(c_{l+1}-1\) are either \(C_*-1\) or \(C_*\). Moreover, Lemma~\ref{lem:reduced-prop} says that at least one of \(\eta_r,\theta_r\) is equal to \(1\), while the weights \((m-r)/m\) and \(r/m\) both lie in \([1/m,(m-1)/m]\). The required one-step comparison therefore reduces to the elementary coefficient checks
\[
z\preceq q,\qquad
\lambda z+(1-\lambda)\preceq q\quad(\lambda\ge1/m),\qquad
\lambda+(1-\lambda)zq\preceq q\quad(\lambda\le(m-1)/m),
\]
and their symmetric versions. This proves \eqref{eq:Fl-domination} by induction.

Choose \(0<\varepsilon<1/(2m)\). Since \(q^{c_l-1}\) is the generating polynomial of \(\operatorname{Bin}(c_l-1,1/m)\), \eqref{eq:Fl-domination} gives
\[
\sum_{N\le\varepsilon c_l}[z^N]F_l(z)
\le
\mathbb P\{\operatorname{Bin}(c_l-1,1/m)\le \varepsilon c_l\}
\le C e^{-\kappa c_l}
\le Cc_l^{-1/2},
\]
after enlarging \(C\) to cover bounded \(c_l\).

Finally split the convex combination for \(D_l\) into terms with \(N(\omega)\le\varepsilon c_l\) and terms with \(N(\omega)>\varepsilon c_l\). The first part has total coefficient \(O(c_l^{-1/2})\) and uniformly bounded variation, while the second part has variation at most \(C(\varepsilon c_l)^{-1/2}\) by \eqref{eq:variation-many-factors}. This proves the desired bound.
\end{proof}
\begin{lemma}
 Suppose $T$ is a restrictive tight map. Then,
 \begin{equation*}
 \norm{D_{l+1}-D_l}_1 \le \dfrac{C}{c_l^{1/2}}.
 \end{equation*}
 \label{lem:shift-L1-bound}
\end{lemma}

\begin{proof}
Let \(\Delta_n=\|D_{n+1}-D_n\|_1\). For \(0\le r\le m\), set \(M_{l,r}=\frac{m-r}{m}D_l+\frac r mD_{l+1}\), with \(M_{l,m}=D_{l+1}\). The recursion also holds for \(r=m\), and Lemmas~\ref{lem:conv-bound}, \ref{lem:tot-var-bound}, and~\ref{lem:c-diff-bound} give
\[
\|D_{ml+r}-M_{l,r}\|_1\le C_1c_l^{-1/2}\qquad(0\le r\le m).
\]
If \(n=ml+r\), then
\[
\Delta_n\le \|D_{ml+r+1}-M_{l,r+1}\|_1+\|D_{ml+r}-M_{l,r}\|_1+\frac1m\Delta_l
\le 2C_1c_l^{-1/2}+\frac1m\Delta_l.
\]
A strong induction, using \(c_{ml+r}\le c_l+1\le2c_l\) and \(m\ge3\), closes the estimate after choosing \(C\) sufficiently large.
\end{proof}

We conclude with some quantitative estimates for $b_l$ and $c_l$.

\begin{lemma}
\label{lem:bl-cl-prop}
 \begin{gather}
 \label{eq:b-diff-1} |b_l - b_{l+1}| \le 1 \\
 \label{eq:c-diff-1} |c_l - c_{l+1}| \le 1 \\
 \label{eq:b-c-comp} b_l \le 2c_l - 1 \\
 \label{eq:c-b-comp} c_l \le b_l
 \end{gather}
\end{lemma}

\begin{proof}
The difference bound for \(b_l\) follows by induction from the recursion: the only nonconstant transitions are the endpoint transitions and the remaining ones are either equal or controlled by the induction hypothesis. The difference bound for \(c_l\) is Lemma~\ref{lem:c-diff-bound}, and \(c_l\le b_l\) is immediate from the definitions. Finally, \(b_l\le2c_l-1\) follows by induction: the case \(r=0\) is inherited from \(l\), and for \(r>0\),
\[
b_{ml+r}\le \max(b_l,b_{l+1})+1\le \min(b_l,b_{l+1})+2
\le 2\min(c_l,c_{l+1})+1=2c_{ml+r}-1.
\]
\end{proof}

\begin{lemma}
We have $c_l\leq \log_ml+2$ whenever $l \ge 2$. This implies $b_l \le 2 \log_m (l+1) + 3$ by \eqref{eq:b-c-comp}.
 \label{lem:bl-size}
\end{lemma}

\begin{proof}
 This holds for $1 \le l \le m$ since $c_l \le 2$ in this range. Assuming $c_l \le \log_ml + 2$, we have
 \[ 
 c_{ml + r} \le c_l + 1 \le \log_{m}l + 3 = \log_m(ml) + 2\le \log_m(ml+r)+2.
 \]
 Thus our claim follows by induction on $l$.
\end{proof}

\section{Construction of Exceptional Set: Proof of Main Theorems}
\label{sec:main-proof}

\subsection{Construction of Exceptional Sets}
\begin{lemma}
 Let $T$ be a restrictive tight map. Recall that $P_n=\set{l\in \mb{N}\,:\, d_l(n)>0}$. If $P_n\ne\emptyset$, then there exists some $l_k = l_k(n) \in P_n$ such that
 \[ 
 P_n \subseteq [l_k -Cb_{l_k} /a_k , l_k + Cb_{l_k}/a_k].
 \]
 Further, $b_{l_k} \le C \log(n+2)$.
 \label{lem:Pn-bound}
\end{lemma}

\begin{proof}
 For any $l \in P_n$, we have
 \[ 
 D_l (n - a_k l) \ne 0 \implies -Rb_l \le n - a_k l \le Rb_l
 \]
 by Lemma~\ref{lem:supp-bound}. In other words,
 \[ 
 a_k l - Rb_l \le n \le a_k l + Rb_l.
 \]
 Now choose $l_k \in P_n$ such that $b_{l_k}$ is the largest possible. Since $|n - a_kl_k| \le Rb_{l_k}$, we can see that $b_{l_k} \le C \log (n+2)$ by Lemma~\ref{lem:bl-size}. Therefore, for any $l \in P_n$, we have
 \[ 
 a_kl_k - Cb_{l_k} \le n - Rb_{l_k} \le a_kl \le n + Rb_l \le n + Cb_{l_k} \le a_k l_k + Cb_{l_k}.
 \]
\end{proof}

For \(l,q\in\mathbb N\), set
\begin{equation}
\label{eq:FG-envelope-def}
F(l,q):=\min_{\substack{-1\le j\le1\\ t\in\mathbb Z,\ |t|\le q}}
D_{l+j}\!\left(\cdot-\frac{t}{m-1}\right),
\qquad
G(l,q):=\max_{\substack{-1\le j\le1\\ t\in\mathbb Z,\ |t|\le q}}
D_{l+j}\!\left(\cdot-\frac{t}{m-1}\right),
\end{equation}
where terms with negative indices are omitted.

\begin{lemma}[Iterated envelope]
\label{lem:iterated-envelope}
Let \(T\) be a restrictive tight map. There is a constant \(C\) such that, for every \(p\ge0\), \(L\ge0\), and integer \(r\) with
\(|r|\le m^p\) and \(m^pL+r\ge0\), if
\[
Q=\lceil Rm(p+1)\rceil,
\]
then
\begin{equation}
\label{eq:iterated-envelope-bound}
F(L,Q)\le D_{m^pL+r}\le G(L,Q).
\end{equation}
Moreover,
\begin{equation}
\label{eq:iterated-envelope-bc}
b_L\le b_{m^pL+r}+C(p+1),
\qquad
c_L\ge c_{m^pL+r}-C(p+1).
\end{equation}
\end{lemma}
\begin{proof}
Put \(N=m^pL+r\). Each use of \eqref{eq:recursion} replaces an index of the form \(ma+s\), \(0\le s<m\), by either \(a\) or \(a+1\). Hence, after \(u\) backward steps, every possible ancestor \(i_u\) satisfies
\[
|m^u i_u-N|\le m^u-1.
\]
Indeed, the error is obtained by adding at most one carry at each of the \(u\) base-\(m\) positions. For \(u=p\), the assumptions \(N=m^pL+r\) and \(|r|\le m^p\) imply
\[
-2m^p+1<m^p(i_p-L)<2m^p,
\]
so \(i_p\in\{L-1,L,L+1\}\). If one of these indices is negative, the corresponding term is simply omitted in the definition of \(F(L,Q)\) and \(G(L,Q)\).

At each backward step the convolution shift is an atom of some \(\alpha_s\) or \(\beta_s\). By the definition of \(R\), the accumulated shift after \(p\) steps has the form \(t/(m-1)\) with
\[
|t|\le R(m-1)p\le Q.
\]
Thus \(D_N\) is a convex combination of functions appearing in the finite family used to define \(F(L,Q)\) and \(G(L,Q)\), which proves \eqref{eq:iterated-envelope-bound}.

Finally, one backward step changes the relevant \(b\)-index by at most an additive constant and changes the relevant \(c\)-index by at most an additive constant, by \eqref{eq:bn-def}, \eqref{eq:cn-def}, Lemma~\ref{lem:c-diff-bound}, and the analogous difference bound for \(b_l\) in Lemma~\ref{lem:bl-cl-prop}. Iterating for \(p\) steps gives \eqref{eq:iterated-envelope-bc}, after enlarging \(C\) to cover \(p=0\).
\end{proof}
\begin{lemma}
\label{lem:FG-bound}
Let \(T\) be a restrictive tight map and assume \(P_n\ne\emptyset\). Then, for all sufficiently large $n$, there exist integers
\(p=p(n)\ge0\), \(L=L(n)\ge0\), and \(Q=Q(n)\ge0\), with
\(p\le C\log(b_{l_k}+2)\) and \(Q=\lceil Rm(p+1)\rceil\), such that
\[
F_n:=F(L,Q),\qquad G_n:=G(L,Q)
\]
satisfy, with the convention \(D_l\equiv0\) for \(l<0\),
\[
F_n(n-a_kl)\le D_l(n-a_kl)\le G_n(n-a_kl)
\qquad(l\in\mathbb Z).
\]
Moreover, the construction gives \(c_L\ge c_{l_k}-C(p+1)\).
\end{lemma}
\begin{proof}
Choose
\[
p=\max\left(0,\left\lceil\log_m\frac{C_0(b_{l_k}+1)}{a_k}\right\rceil\right),
\qquad
L=\left\lfloor \frac{l_k}{m^p}+\frac12\right\rfloor,
\]
with \(C_0\) large, and set \(Q=\lceil Rm(p+1)\rceil\). The rounding gives \(|l_k-m^pL|\le m^p/2\). By Lemma~\ref{lem:Pn-bound}, and by increasing \(C_0\), every \(l\in P_n\) satisfies \(|l-m^pL|\le m^p\). Lemma~\ref{lem:iterated-envelope} therefore gives
\[
D_l(n-a_kl)\le G(L,Q)(n-a_kl)\qquad(l\in P_n),
\]
while the upper bound is trivial for \(l\notin P_n\).

For the lower bound, suppose first that \(l\ge0\) and \(F(L,Q)(n-a_kl)>0\). Since the unshifted term \(D_L\) occurs in the defining minimum for \(F(L,Q)\), Lemma~\ref{lem:supp-bound} gives
\[
|n-a_kl|\le Rb_L.
\]
Also \(l_k\in P_n\), so \(|n-a_kl_k|\le Rb_{l_k}\). Hence
\[
a_k|l-l_k|\le R(b_L+b_{l_k}).
\]
Applying \eqref{eq:iterated-envelope-bc} to \(l_k=m^pL+(l_k-m^pL)\) gives \(b_L\le b_{l_k}+C(p+1)\). The choice of \(p\) and a sufficiently large \(C_0\) then imply
\[
R(b_L+b_{l_k})\le \frac12 a_km^p.
\]
Together with \(|l_k-m^pL|\le m^p/2\), this yields \(|l-m^pL|\le m^p\). Lemma~\ref{lem:iterated-envelope} now gives
\[
F(L,Q)(n-a_kl)\le D_l(n-a_kl).
\]
If \(l<0\), then \(D_l\equiv0\). Since the unshifted term \(D_L\) occurs in the defining minimum for \(F(L,Q)\),
\[
\supp F(L,Q)\subset \supp D_L\subset[-Rb_L,Rb_L]\subset[-C\log(n+2),C\log(n+2)].
\]
As \(n-a_kl\ge n\), the lower inequality is trivial for all sufficiently large \(n\). Finally, applying \eqref{eq:iterated-envelope-bc} with \(r=l_k-m^pL\) gives \(c_L\ge c_{l_k}-C(p+1)\), and the bound \(p\le C\log(b_{l_k}+2)\) follows from the definition of \(p\).
\end{proof}
Thus, by the methods discussed in Section~\ref{sec:exceptional-set}, we can find an exceptional set for $(A_k, A_k)$. We proceed to show an upper bound of its size by proving that only a small number of $F_n, G_n$ behave badly. We do this by showing that $c_l$ being large implies good conditions on $(F_n, G_n)$ (Lemma~\ref{lem:FG-good}) and then showing that only a small number of $c_l$ can be small (Lemma~\ref{lem:small-cl-bound}).

\begin{lemma}
 Let $T$ be a restrictive tight map. Then for any $l\ge0$ and $q\ge1$,
 \begin{gather}
 \label{eq:FG-L1} \|G(l, q) - F(l, q)\|_1 \le \frac{C q^2}{c_l^{1/2}} \\
 \label{eq:F-var}V(F(l, q)) \le \frac{Cq}{c_l^{1/2}} \\
 \label{eq:G-var}V(G(l, q)) \le \frac{Cq}{c_l^{1/2}} .
 \end{gather}
 \label{lem:FG-good}
\end{lemma}

\begin{proof}
For each non-omitted pair \(-1\le j\le1\) and \(t\in\mathbb Z\) with \(|t|\le q\), Lemma~\ref{lem:shift-L1-bound}, Lemma~\ref{lem:tot-var-bound}, Lemma~\ref{lem:conv-bound}, and \(|c_{l+j}-c_l|\le C\) give
\[
\left\|D_{l+j}\left(\cdot-\frac{t}{m-1}\right)-D_l\right\|_1
\le \frac{C(|t|+1)}{c_l^{1/2}}.
\]
Taking finite maxima and minima and summing these bounds yields
\[
\|G(l,q)-F(l,q)\|_1
\le \frac{C}{c_l^{1/2}}
\sum_{\substack{t\in\mathbb Z\\ |t|\le q}}(|t|+1)
\le \frac{Cq^2}{c_l^{1/2}}.
\]
For the variation estimates, the Appendix bound for finite maxima/minima and Lemma~\ref{lem:tot-var-bound} give
\[
V(G(l,q)),\ V(F(l,q))
\le
\sum_{\substack{-1\le j\le1,\ l+j\ge0\\ t\in\mathbb Z,\ |t|\le q}}
V\left(D_{l+j}\left(\cdot-\frac{t}{m-1}\right)\right)
\le \frac{Cq}{c_l^{1/2}}.
\]
\end{proof}
\begin{corollary}
\label{cor:Jh-exceptional}
 Let $T$ be a restrictive tight map and take any increasing function $h:\mb{R}^+\to \mb{R}^+$ diverging to infinity. Then,
 \[
 J(h) = \{ n \in \mathbb{N}: P_n=\emptyset \}\cup\{ n \in \mathbb{N}: P_n\ne\emptyset,\ c_{l_k(n)} \le h(n) \}
 \]
 is an exceptional set for $(A_k, A_k)$.
\end{corollary}

\begin{proof}
For \(n\notin J(h)\), Lemma~\ref{lem:FG-bound} supplies \(F_n\le D_l\le G_n\) on the relevant lattice points. Moreover \(Q=O(p+1)\),
\[
p+1=O(\log(c_{l_k}+2)),\qquad
c_L\ge c_{l_k}-C\log(c_{l_k}+2).
\]
Also \(0\le F_n\le D_L\le G_n\) and \(\|D_L\|_1=1\), hence \(\|F_n\|_1\le1\le\|G_n\|_1\).
Since \(c_{l_k}>h(n)\to\infty\), Lemma~\ref{lem:FG-good} gives \(\|G_n-F_n\|_1,V(F_n),V(G_n)\to0\). The sandwich estimate from Proposition~\ref{prop:Fn-Gn} then gives \(\mu(A_k\cap T^{-n}A_k)\to\mu(A_k)^2\).
\end{proof}

\subsection{\texorpdfstring{Upper bounds for \(J(h)\)}{Upper bounds for J(h)}}
\label{subsec:Jh-size}

We proceed to give an upper bound of $|J(h) \cap [0, n]|$ for restrictive tight maps. This proves Theorem~\ref{thm:general}, and gives Corollary~\ref{cor:chacon} as a corollary.

\begin{lemma}
Let \(N_{M, q} := | \{ m^{q-1} \le l < m^{q}: c_l \le M\}|\). Then, for some constant \(C_m>0\) depending only on \(m\),
\[
N_{M,q}\le (C_mq)^M.
\]
Moreover, for \(M\ge2\),
\[
N_{M,q}\ge \binom{q-1}{M-2}.
\]
\label{lem:small-cl-bound}
\end{lemma}
\begin{proof}
Writing \(L=ml+r\), with \(0\le r<m\), gives
\[
c_{ml}=c_l,
\qquad
c_{ml+r}=\min(c_l,c_{l+1})+1\quad(1\le r<m).
\]
If \(m^q\le L<m^{q+1}\), then \(m^{q-1}\le l<m^q\). The terms with \(r=0\) contribute at most \(N_{M,q}\). For \(1\le r<m\), the condition \(c_{ml+r}\le M\) implies
\[
c_l\le M-1\quad\text{or}\quad c_{l+1}\le M-1.
\]
Thus each bad index at level \(q\) can be charged to either \(l\) or \(l+1\); the only index not lying in the interval \([m^{q-1},m^q)\) is the endpoint \(l+1=m^q\). Hence
\[
N_{M,q+1}\le N_{M,q}+2(m-1)(N_{M-1,q}+1).
\]
Starting from \(N_{0,q}=0\) and the bounded initial values \(N_{M,1}=O_m(1)\), a double induction on \(M\) and \(q\) gives \(N_{M,q}\le(C_mq)^M\) after increasing \(C_m\).

For the lower bound, the \(r=0\) branch gives \(N_{M,q}\) inherited terms. In addition, if \(c_l\le M-1\), then for every \(1\le r<m\),
\[
c_{ml+r}=\min(c_l,c_{l+1})+1\le M.
\]
Therefore
\[
N_{M,q+1}\ge N_{M,q}+(m-1)N_{M-1,q}\ge N_{M,q}+N_{M-1,q}.
\]
The bound \(N_{M,q}\ge {q-1\choose M-2}\) follows from this recurrence, the base case \(N_{2,q}\ge1\), and Pascal's identity.
\end{proof}
\begin{lemma}
\label{lem:empty-times-small-c}
Let $T$ be a restrictive tight map and fix $k$. There exists a constant $C_k$ such that the following holds. If $P_n=\emptyset$ and $\ell=\ell(n)$ is any nonnegative integer with
\begin{equation}
\label{eq:nearest-ell}
|n-a_k\ell|\le a_k/2+1,
\end{equation}
then $c_\ell\le C_k$.
\end{lemma}

\begin{proof}
Constants may depend on \(k\). Choose
\[
p=\max\left(1,\left\lceil\log_m\frac{C_0(b_\ell+2)}{a_k}\right\rceil\right),
\quad
L=\left\lfloor \frac{\ell}{m^p}+\frac12\right\rfloor,
\quad
Q=\lceil Rm(p+1)\rceil,
\]
with \(C_0\) large, and put \(F=F(L,Q)\). Lemma~\ref{lem:iterated-envelope} gives \(F\le D_j\) whenever \(|j-m^pL|\le m^p\), and
\[
b_L\le b_\ell+C(p+1),\qquad c_L\ge c_\ell-C(p+1).
\]
The choice of \(p\) gives \(Rb_L+a_k/2+1\le a_km^p/2\). Since the unshifted term \(D_L\) occurs in the defining minimum of \(F\),
\[
\supp F\subset \supp D_L\subset[-Rb_L,Rb_L]\subset[-C\log(n+2),C\log(n+2)],
\]
so \(F(n-a_kj)=0\) for \(j<0\) and all sufficiently large \(n\). Hence, if \(F(n-a_kj)>0\), then \(j\ge0\), and the support bound and \eqref{eq:nearest-ell} imply \(|j-m^pL|\le m^p\), so \(F(n-a_kj)\le D_j(n-a_kj)\). Since \(P_n=\emptyset\), \(\sum_{j\in\mathbb Z}F(n-a_kj)=0\). The Appendix summation estimate yields \(a_k^{-1}\|F\|_1\le V(F)\). Lemma~\ref{lem:FG-good} gives
\[
V(F)\le C(p+1)c_L^{-1/2},
\qquad
\|D_L-F\|_1\le C(p+1)^2c_L^{-1/2}.
\]
Since \(\|D_L\|_1=1\), either directly or through \(\|F\|_1\ge1/2\) we get \(c_L\le C_k(p+1)^4\). Combining this with \(c_L\ge c_\ell-C(p+1)\) and \(p+1\le C_k\log(c_\ell+2)\), we obtain \(c_\ell\le C_k\log^4(c_\ell+2)\), which forces \(c_\ell\) to be bounded. Enlarging \(C_k\) handles the finitely many small \(n\).
\end{proof}

\begin{lemma}
 Suppose $T$ is a restrictive tight map. For any increasing function $h:\mb{R}^+\to \mb{R}^+$ diverging to infinity, we have, for all sufficiently large $M$,
 \begin{equation}
 |J(h) \cap [0, M]| \le C h(M)(C\log(M+2))^{\lceil h(M)\rceil+1}.
 \end{equation}
 \label{lem:Jh-bound}
\end{lemma}

\begin{proof}
Assign each \(n\in J(h)\cap[0,M]\) to an index \(l\): use \(l=l_k(n)\) if \(P_n\ne\emptyset\), and otherwise choose \(l\) with \(|n-a_kl|\le a_k/2+1\). Lemma~\ref{lem:supp-bound} and Lemma~\ref{lem:empty-times-small-c} show that, for large \(M\), every assigned index satisfies \(c_l\le h(M)\), and each fixed index receives only \(O(h(M))\) values of \(n\). Also \(n\le M\), the support bound, and \(b_l\le C\log(l+2)\) imply \(l\le C(M+1)/a_k\). Therefore
\[
\begin{aligned}
|J(h)\cap[0,M]|
&\le Ch(M)\left|\left\{l\le C(M+1)/a_k:\ c_l\le h(M)\right\}\right| \\
&\le Ch(M)\sum_{1\le q\le C\log(M+2)}N_{\lceil h(M)\rceil,q} \\
&\le Ch(M)(C\log(M+2))^{\lceil h(M)\rceil+1}.
\end{aligned}
\]
\end{proof}

\begin{corollary}
\label{cor:J_k}
 Suppose $T$ is a restrictive tight map. For any increasing $h : \mathbb{R}^+ \to \mathbb{R}^+$ diverging to infinity, we may choose an exceptional set $J_k$ of $(A_k, A_k)$ such that
 \[
 |J_k \cap [0, n]| \le (\log n)^{h(n)}
 \]
 for all sufficiently large $n$.
\end{corollary}

\begin{proof}
 Recall Corollary~\ref{cor:Jh-exceptional} and Lemma~\ref{lem:Jh-bound}. By choosing an appropriate $\bar{h}$ according to $h$, we can ensure that
 \[ 
 C\bar{h}(n)(C\log(n+2))^{\lceil\bar{h}(n)\rceil+1}
 \le (\log n)^{h(n)}.
 \]
 Thus setting $J_k = J(\bar{h})$ completes the proof.
\end{proof}

\begin{proof}[Proof of Theorem~\ref{thm:general}]
Let
\[
\mcr{C}:=\{T^rA_k:\ k\in\mathbb N,\ 0\le r<h_k\},
\]
where \(h_k\) is the height of the \(k\)-th tower. For \(T^iA_k,T^jA_\ell\in\mcr C\), choose \(K\) so that both are finite disjoint unions of levels of the \(K\)-th tower:
\[
T^iA_k=\bigsqcup_{u\in U}T^uA_K,
\qquad
T^jA_\ell=\bigsqcup_{v\in V}T^vA_K,
\]
where \(U,V\subset\{0,1,\ldots,h_K-1\}\), modulo null sets. Since
\[
\mu(T^uA_K\cap T^{-n}T^vA_K)=\mu(A_K\cap T^{-(n+u-v)}A_K),
\]
finite unions of shifted exceptional sets for \((A_K,A_K)\) give exceptional sets for every pair in \(\mcr C\). Shifts and finite unions only change the counting estimate by constants and \(O(1)\). For each resulting pair, these constants are absorbed after discarding a finite initial segment of the associated exceptional set and applying Corollary~\ref{cor:J_k} with a slightly smaller divergent function.
Indeed, choose divergent \(h_0,h_1\) so slowly that \(h_1(n)(\log n)^{h_0(n)}\le(\log n)^{h(n)}\) eventually; use \(h_0\) for the pairwise bounds and \(h_1\) in Corollary~\ref{prop:J_C}. Corollary~\ref{prop:J_C} therefore yields an exceptional set \(J\) for \(\mcr C\) with
\[
|J\cap[0,n]|\le(\log n)^{h(n)}
\]
for all sufficiently large \(n\).

Finally, by the standing rank-one convention, \(\mcr C\) generates \(\mcr B\) modulo null sets. More explicitly, the finite level partitions
\[
\mathcal P_k:=\{T^rA_k:0\le r<h_k\}
\]
refine along the construction, their tower unions exhaust \(X\) modulo null sets, and finite disjoint unions of their atoms approximate every set in \(\mcr B\) in measure. These atoms are members of \(\mcr C\), so the finite disjoint approximation hypothesis of Proposition~\ref{prop:J_B} is satisfied. Hence \(J\) is exceptional for \(\mcr B\).
\end{proof}

By Proposition~\ref{prop:J_f}, this leads to the following corollary.
\begin{corollary}
 For restrictive tight maps, the exceptional set $J$ in Theorem~\ref{thm:general} is exceptional for every $f,g\in L^2(\mu)$. 
\end{corollary}

\subsection{Lower bound on Exceptional Set}

In this section, we prove Theorem~\ref{thm:lower-bound} to show that the upper bound in Theorem~\ref{thm:general} is optimal in some sense. Throughout this subsection, $T$ is a tight map whose spacer sequence satisfies $s_{m-1}=0$. Recall the definition of $R$ in \eqref{eq:R-def} and the definition $a_k = \mu(A_k)^{-1}$.

We use the recursively defined sequences $b_l,c_l$ from \eqref{eq:bn-def}--\eqref{eq:cn-def}. The lower-bound argument below uses the support estimate Lemma~\ref{lem:supp-bound} and the purely combinatorial estimates for $b_l,c_l$; it does not use the restrictive hypotheses except through the stated assumption $s_{m-1}=0$.

\begin{lemma}
Assume in addition that $m\ge3$, and let
\[
E_k=\{n\in\mathbb N:\ \mu(A_k\cap T^{-n}A_k)=0\}.
\]
For \(a_k\) sufficiently large,
\[
|E_k\cap[0,a_kn]|
\ge
\binom{\lfloor \log_m n \rfloor-1}{\lfloor a_k/(8R)\rfloor-2}.
\]
\label{lem:lower-bound}
\end{lemma}

\begin{proof}
By Lemma~\ref{lem:supp-bound}, \(d_l\) is supported in
\[
I_l:=[a_kl-Rb_l,\ a_kl+Rb_l].
\]
For \(a_k>R\), the inequalities \(|b_l-b_{l+1}|\le1\) imply
\[
a_k(l+1)-Rb_{l+1}\ge a_kl-Rb_l,
\qquad
a_k(l+1)+Rb_{l+1}\ge a_kl+Rb_l,
\]
so the endpoints of \(I_l\) are increasing in \(l\). If \(c_l\le\lfloor a_k/(8R)\rfloor\), then \(b_l\le2c_l-1\le a_k/(4R)\), and \(b_{l+1}\le b_l+1\). Hence the gap between \(I_l\) and \(I_{l+1}\) has length at least
\[
a_k-R(b_l+b_{l+1})\ge \frac{a_k}{2}-R>1
\]
for \(a_k\) large enough. Choose an integer \(n_l\) in this gap. Then \(n_l\notin I_j\) for every \(j\), so \(d_j(n_l)=0\) for all \(j\), and hence \(n_l\in E_k\). Distinct \(l\)'s give distinct gaps. For the chosen integer in the gap, \(n_l<a_k(l+1)-Rb_{l+1}<a_k(l+1)\). Thus \(l<n\) implies \(l+1\le n\), hence \(n_l<a_kn\). Therefore
\[
|E_k\cap[0,a_kn]|
\ge
\left|\{1\le l<n:\ c_l\le \lfloor a_k/(8R)\rfloor\}\right|.
\]
The desired lower bound follows from Lemma~\ref{lem:small-cl-bound}.
\end{proof}

\begin{proof}[Proof of Theorem~\ref{thm:lower-bound}]
First suppose $m=2$. Since $s_1=0$, outside the null set of points eventually remaining in the last subcolumn, each return from $A_k$ crosses one $k$-tower and the $s_0$ spacers, so the first return time is $r_k=h_k+s_0=a_k$. Hence $t_l=la_k$ a.e. for every $l\ge0$. Therefore
\[
E_k:=\{n\in\mathbb N:\ \mu(A_k\cap T^{-n}A_k)=0\}
\]
contains all integers which are not of the form $a_kl$ with $l\ge0$. Choose $k$ with $a_k\ge2$ and set $A=B=A_k$. Then $|E_k\cap[0,N]|\ge N/2-O(1)$. If $J_{A,A}$ is exceptional for $(A,A)$, then Lemma~\ref{lem:finite-modification-exceptional}(2), applied with $\tau=\mu(A_k)^2/2$, gives that $E_k\setminus J_{A,A}$ is finite. The desired lower bound follows.

Now assume $m\ge3$. Choose \(k\) so large that
\[
\left\lfloor\frac{a_k}{8R}\right\rfloor-2>t+1.
\]
Let \(A=B=A_k\). By Lemma~\ref{lem:lower-bound}, and by taking \(n=\lfloor N/a_k\rfloor\), we obtain
\[
|E_k\cap[0,N]|\ge C(\log N)^{t+1}
\]
for all sufficiently large \(N\), after changing \(C>0\). If \(J_{A,A}\) is exceptional for \((A,A)\), then Lemma~\ref{lem:finite-modification-exceptional}(2), applied with \(\tau=\mu(A_k)^2/2\), gives that \(E_k\setminus J_{A,A}\) is finite. Hence
\[
|J_{A,A}\cap[0,N]|\ge (\log N)^t
\]
for all sufficiently large \(N\). This implies the stated bound.
\end{proof}

\section{Applications and Related Problems}
\label{sec:applications}

\subsection{\texorpdfstring{Generalization to \(\mathbb{R}\) and \(\mathbb{Z}^d\)--actions}{Generalization to R and Zd-actions}}
\label{subsec:gen-action}

While weak mixing is most commonly defined for $\mathbb{Z}$--actions, the definition extends naturally to flows and higher-rank actions. For a measure-preserving flow $(T^t)_{t\in\mathbb{R}}$ on $(X,\mcr{B},\mu)$, one sets
\[
C_T \;:=\;\frac{1}{T}\int_{0}^{T}\bigl|\mu(A\cap T^{-t}B)-\mu(A)\,\mu(B)\bigr|\,dt
\;\longrightarrow\;0
\quad (T\to\infty)
\]
for all measurable $A,B\subseteq X$. The vanishing of these Cesàro averages for every $A,B$ characterizes weak mixing of the flow.

Similarly, a measure-preserving $\mathbb{Z}^d$--action $(T^n)_{n\in\mathbb{Z}^d}$ is weak mixing if
\[
\frac{1}{(2N+1)^d}\sum_{n\in[-N,N]^d}\bigl|\mu(A\cap T^{-n}B)-\mu(A)\,\mu(B)\bigr|
\;\longrightarrow\;0
\quad (N\to\infty),
\]
for every pair of measurable sets $A,B\in\mcr{B}$.

In each case, the failure of mixing is measured by exceptional sets. For a discrete-time transformation one obtains a zero-density exceptional set $J_{A,B}\subseteq\mathbb{N}$. In the continuous-time setting the corresponding exceptional subset of $[0,\infty)$ has Lebesgue measure $o(T)$ in $[0,T]$. For a $\mathbb{Z}^d$--action one obtains $J_{A,B}\subseteq\mathbb{Z}^d$ whose proportion in the cubes $[-N,N]^d$ vanishes as $N\to\infty$.

Theorem \ref{prop:weak-conv-rate2} is the $\mb{R}$ and $\mb{Z}^d$--action version of Proposition \ref{prop:weak-conv-rate}.

\begin{theorem}
\label{prop:weak-conv-rate2}
Fix $p\ge1$ and let $b_T\to 0$ be a positive sequence (or function).

\begin{enumerate}
 \item\textbf{($\mb{R}$-\emph{action})} 
 Let $(X,\mcr{B},\mu,(\phi^{s})_{s\in\mathbb{R}})$ be a measure-preserving flow and let $A,B\in\mcr{B}$ satisfy
 \[
 \frac{1}{T}\int_{0}^{T}\bigl|\mu\bigl(A\cap (\phi^{s})^{-1}B\bigr)-\mu(A)\mu(B)\bigr|^{p}\,ds
 =o\bigl(b_T\bigr)
 \quad(T\to\infty).
 \]
 Then there exists an exceptional set $J_{A,B}\subseteq[0,\infty)$ such that
 \[
 \operatorname{Leb}\bigl(J_{A,B}\cap[0,T]\bigr)
 =o\bigl(T\,b_T\bigr)
 \quad\text{and}\quad
 \mu\bigl(A\cap (\phi^{s})^{-1}B\bigr)\longrightarrow\mu(A)\mu(B)
 \]
 as $s\to\infty$ with $s\notin J_{A,B}$. 

 \item\textbf{$(\mathbb Z^d$--\emph{action})}
 Let $(X,\mcr{B},\mu,(T^n)_{n\in\mathbb Z^d})$ be a measure-preserving $\mathbb Z^d$--action and let $A,B\in\mcr{B}$ satisfy
 \[
 \frac{1}{(2N+1)^d}\sum_{n\in[-N,N]^d}\bigl|\mu(A\cap T^{-n}B)-\mu(A)\mu(B)\bigr|^{p}
 =o\bigl(b_N\bigr)
 \quad(N\to\infty).
 \]
 Then there exists an exceptional set $J_{A,B}\subseteq\mathbb Z^d$ such that
 \[
 \bigl|J_{A,B}\cap[-N,N]^d\bigr|
 =o\bigl((2N+1)^d b_N\bigr)
 \quad\text{and}\quad
 \mu(A\cap T^{-n}B)\longrightarrow\mu(A)\mu(B)
 \]
 as $|n|\to\infty$ with $n\notin J_{A,B}$.
\end{enumerate}
\end{theorem}

\begin{proof}[Proof of Theorem~\ref{prop:weak-conv-rate2}]
Throughout, put
\[
 a(s)=\bigl|\mu\bigl(A\cap(\phi^{s})^{-1}B\bigr)-\mu(A)\mu(B)\bigr|^{p}\quad(s\ge0),\qquad
 a_{n}=\bigl|\mu(A\cap T^{-n}B)-\mu(A)\mu(B)\bigr|^{p}.
\]

\medskip
\noindent\textbf{(1) \(\mathbb R\)-action.}
For \(k\in\mathbb N\), set
\[
E_k:=\{s\ge0:\ a(s)>1/k\}.
\]
Markov's inequality gives
\[
\operatorname{Leb}(E_k\cap[0,T])
\le
k\int_0^T a(s)\,ds
=
o(Tb_T).
\]
Choose \(T_k\to\infty\) so that, for all \(T\ge T_k\),
\[
\operatorname{Leb}(E_k\cap[0,T])\le \frac1k Tb_T.
\]
Define
\[
J_{A,B}:=\bigcup_{k=1}^{\infty}E_k\cap[T_k,T_{k+1}).
\]
If \(T\in[T_K,T_{K+1})\), then
\[
J_{A,B}\cap[0,T]\subseteq E_K\cap[0,T],
\]
so
\[
\operatorname{Leb}(J_{A,B}\cap[0,T])\le \frac1K Tb_T=o(Tb_T).
\]
Moreover, if \(s\notin J_{A,B}\) and \(s\in[T_K,T_{K+1})\), then \(a(s)\le1/K\). Hence \(a(s)\to0\) as \(s\to\infty\) outside \(J_{A,B}\).

\medskip
\noindent\textbf{(2) \(\mathbb Z^d\)-action.}
For \(k\in\mathbb N\), set
\[
E_k:=\{n\in\mathbb Z^d:\ a_n>1/k\}.
\]
Then
\[
|E_k\cap[-N,N]^d|
\le
k\sum_{n\in[-N,N]^d}a_n
=
o((2N+1)^db_N).
\]
Choose \(N_k\to\infty\) such that, for all \(N\ge N_k\),
\[
|E_k\cap[-N,N]^d|
\le
\frac1k(2N+1)^db_N.
\]
Define
\[
J_{A,B}:=
\bigcup_{k=1}^{\infty}
E_k\cap\bigl([-N_{k+1},N_{k+1}]^d\setminus[-N_k,N_k]^d\bigr).
\]
If \(N\in[N_K,N_{K+1})\), then
\[
J_{A,B}\cap[-N,N]^d\subseteq E_K\cap[-N,N]^d,
\]
and hence
\[
|J_{A,B}\cap[-N,N]^d|
\le
\frac1K(2N+1)^db_N=o((2N+1)^db_N).
\]
Finally, outside \(J_{A,B}\) we have \(a_n\le1/K\) on the \(K\)-th annulus, so \(a_n\to0\) as \(|n|\to\infty\).
\end{proof}

Hence, we can find an upper bound on the size of the exceptional set given the rate of weak mixing. 

\subsection{\texorpdfstring{Unified Applications of Proposition~\ref{prop:weak-conv-rate}~and~Theorem~\ref{prop:weak-conv-rate2} across weak mixing models}{Unified Applications across weak mixing models}}
\label{subsec:applications}

When a cited quantitative weak-mixing estimate has the form \(O(r_R)\), we apply Proposition~\ref{prop:weak-conv-rate} or Theorem~\ref{prop:weak-conv-rate2} with any \(b_R\) satisfying \(r_R=o(b_R)\); endpoint exponents are therefore replaced by arbitrary smaller ones. The set-valued statements below assume that the centered indicators \(\mathbf1_A-\mu(A)\) and \(\mathbf1_B-\mu(B)\) lie in the regularity class covered by the cited estimate. Equivalently, the same arguments apply to arbitrary zero-mean observables in that class.

In this subsection we illustrate how the discrete-time, continuous-time, and $\mathbb Z^d$ versions of our convergence theorem yield concrete exceptional-set estimates in five paradigmatic settings: random substitution tilings (flows); interval exchange transformations (IETs); translation flows; primitive substitution $\mathbb Z$--actions; and self-affine substitution tilings ($\mathbb Z^d$--actions).

A \emph{random substitution tiling} is the tiling space $(\Omega_x,\phi^t,\mu_x)$ obtained by choosing, according to a shift-invariant ergodic measure $\mu$ on $\{1,\dots,N\}^{\mathbb Z}$, a sequence of compatible uniformly expanding substitutions $S_{x_k}$. Under the hypotheses of Theorem 1.2 in \cite{trevino2020quantitative}, there exists $\alpha'\in(0,1]$ such that for Lipschitz zero-mean observables $f,g$
\[
 \frac1T\int_0^T\!\bigl|\langle f\circ\phi^t, g\rangle\bigr|\,dt=O\!\bigl(T^{-\frac{\alpha'}{2}+\varepsilon}\bigr)\quad(\forall\varepsilon>0).
\]

\begin{corollary}
\label{cor:random-tiling-exceptional}
For sets $A,B\subseteq\Omega_x$ whose centered indicators belong to the regularity class covered by the quoted estimate, and for every $0<\beta<\alpha'/2$, there is an exceptional set $J_{A,B}\subseteq[0,\infty)$ with
\[
 \operatorname{Leb}\bigl(J_{A,B}\cap[0,T]\bigr)=o(T^{1-\beta}).
\]
\end{corollary}

\begin{proof}
Set $f=\mathbf1_A-\mu_x(A)$ and $g=\mathbf1_B-\mu_x(B)$, which belong to this class by assumption. Choose $\varepsilon>0$ so that $\beta<\alpha'/2-\varepsilon$. Then the quoted estimate is $o(T^{-\beta})$, so the continuous-time case of Theorem~\ref{prop:weak-conv-rate2} applies with $p=1$ and $b_T=T^{-\beta}$.
\end{proof}

An \textit{interval exchange transformation} (IET) rearranges subintervals of $[0,1)$ by translations. Avila--Forni--Safaee \cite{IET_rate} prove quantitative weak-mixing estimates for IETs; when these estimates are written with endpoint exponents as $O(N^{-\alpha_0})$ or $O(\log^{-a_0}N)$, we use arbitrary smaller exponents below.

\begin{corollary}
\label{cor:IET}
Let $(I,\mcr B,\mu,T)$ be a typical IET and let $A,B\subseteq I$ be sets whose centered indicators belong to the regularity class required in the cited estimate.
\begin{enumerate}
 \item[\textup{(a)}] If $T$ is non-rotation class and the cited estimate gives $O(N^{-\alpha_0})$ for some $\alpha_0>0$, then for every $0<\alpha<\alpha_0$ one has $|J_{A,B}\cap[0,N]|=o(N^{1-\alpha})$.
 \item[\textup{(b)}] If $T$ is rotation class and the cited estimate gives $O(\log^{-a_0}N)$ for some $a_0>0$, then for every $0<a<a_0$ one has $|J_{A,B}\cap[0,N]|=o\bigl(N\log^{-a}N\bigr)$.
\end{enumerate}
\end{corollary}

\begin{proof}
In case (a), $O(N^{-\alpha_0})=o(N^{-\alpha})$ for every $0<\alpha<\alpha_0$, so take $b_N=N^{-\alpha}$. In case (b), $O(\log^{-a_0}N)=o(\log^{-a}N)$ for every $0<a<a_0$, so take $b_N=\log^{-a}N$. Applying Proposition~\ref{prop:weak-conv-rate} (the $\mathbb Z$--action case) gives the stated bounds.
\end{proof}

Avila--Forni \cite{IET} proved weak mixing for typical translation flows on higher-genus surfaces. Effective weak-mixing and spectral-measure estimates in this setting are available in Forni and Bufetov--Solomyak \cite{Forni2022,BufetovSolomyak2021}; the precise exponent and regularity class depend on the theorem used. The following application is conditional on any quantitative Ces\`aro weak-mixing estimate available for the relevant regularity class.

\begin{corollary}
\label{cor:translation-surface-exceptional}
Let $(X,\phi^t,\mu)$ be the translation flow on a typical genus~$g\ge2$ surface. Suppose that, for the relevant regularity class, a quantitative Ces\`aro weak-mixing estimate gives a bound $O(T^{-\alpha_0})$ for some $\alpha_0>0$. Then for sets $A,B\subseteq X$ whose centered indicators belong to that class and every $0<\alpha<\alpha_0$, there is $J_{A,B}\subseteq[0,\infty)$ with $\operatorname{Leb}(J_{A,B}\cap[0,T])=o(T^{1-\alpha})$.
\end{corollary}

\begin{proof}
Since $O(T^{-\alpha_0})=o(T^{-\alpha})$ for every $0<\alpha<\alpha_0$, apply the continuous-time case of Theorem~\ref{prop:weak-conv-rate2} with $p=1$ and $b_T=T^{-\alpha}$.
\end{proof}

Let $(X_\zeta,T,\mu)$ be the uniquely ergodic system arising from a primitive, aperiodic substitution $\zeta$. Bufetov--Marshall-Maldonado--Solomyak \cite{bufetov2025local} show
\[
 \frac1N\sum_{k=0}^{N-1}\!|\langle U^k f,g\rangle|^2=O\bigl((\log N)^{-\gamma_0}\bigr),
\]
with $\gamma_0>0$ explicit.

\begin{corollary}
\label{cor:substitution-exceptional}
For sets $A,B\subseteq X_\zeta$ whose centered indicators belong to the regularity class required in the cited estimate and every $0<\gamma<\gamma_0$, one has
\[
|J_{A,B}\cap[0,N]|=o\bigl(N(\log N)^{-\gamma}\bigr).
\]
\end{corollary}

\begin{proof}
This is Proposition~\ref{prop:weak-conv-rate} with $p=2$ and $b_N=(\log N)^{-\gamma}$, since the quoted estimate is $o((\log N)^{-\gamma})$ for every $\gamma<\gamma_0$.
\end{proof}

For a self-affine substitution tiling of $\mathbb R^d$, Marshall-Maldonado \cite[Thm.~6.5]{marshall2024quantitative} showed
\[
 \frac{1}{(2N+1)^d}\sum_{n\in[-N,N]^d}|\langle U^n f,g\rangle|^2=O\bigl((\log N)^{-\gamma_0}\bigr).
\]

\begin{corollary}
\label{cor:selfaffine-exceptional-mod}
In the associated $\mathbb Z^d$-action $(X_\zeta,T^n)$, for sets $A,B$ whose centered indicators belong to the regularity class required in the cited estimate and every $0<\gamma<\gamma_0$, there is $J_{A,B}\subseteq\mathbb Z^d$ with
\[
 |J_{A,B}\cap[-N,N]^d|=o\bigl((2N+1)^d(\log N)^{-\gamma}\bigr).
\]
\end{corollary}

\begin{proof}
Apply the $\mathbb Z^d$-action case of Theorem~\ref{prop:weak-conv-rate2} with $p=2$ and $b_N=(\log N)^{-\gamma}$.
\end{proof}

Moll~\cite{moll2023speed}, which references an earlier arXiv version of our draft, proved that for any zero-mean Lipschitz observable \(f\) and \(g\in L^2\),
\[
 \frac{1}{N}\sum_{k=0}^{N-1}\bigl|\langle U^k f,\,g\rangle\bigr|^2
 = O\bigl(\|f\|_L^2\,\|g\|_2^2\,[\log_3 N]^{-1/6}\bigr),
\]
and established a matching lower bound
\[
 \sum_{k=0}^{N-1}\bigl|\langle U^k f_N,\,g_N\rangle\bigr|^2
 \ge C\,\frac{N}{(\log N)^2}\,\|f_N\|_L^{\frac12}\|f_N\|_2^{\frac12}\|g_N\|_2.
\]

In Moll's approach, one then applies Lemma~\ref{lem:weak-conv-speed} (with \(p=2\)) to deduce that, for every \(\gamma<1/6\), his method produces an exceptional set satisfying
\[
 \bigl|J_{f,g}\cap[0,N]\bigr|
 = o(N[\log_3 N]^{-\gamma}).
\]
In contrast, our Theorem~\ref{thm:general} and Corollary~\ref{cor:chacon} give the much stronger bound
\[
 \lvert J\cap[0,n]\rvert \le (\log n)^{h(n)},
\]
highlighting the gap between ``automatic'' exceptional-set estimates derived solely from weak mixing rates and the sharper bounds obtained via direct construction. It therefore seems plausible that, by blending Moll's spectral-measure techniques with our methods, one could further improve exceptional-set bounds in some of the applications discussed above.

\subsection{A polynomial lower bound example in a rank-one transformation}
\label{subsec:rank-one-lower-bound}

We record a self-contained rank-one example with polynomial lower growth of exceptional sets. For $k\geq 0$, let
\[
m_k:=2^{k+1}, \quad 
s_{k,r}:=
\begin{cases}
0,&0\le r<m_k-1,\\
1,&r=m_k-1.
\end{cases}
\]
Let $(X,\mathscr B,\mu,T)$ be the rank-one transformation with cut sequence $\{m_k\}_{k\ge0}$ and spacer sequence $\{s_{k,r}\}$, where $m_k$ denotes the number of subcolumns at stage $k$: at stage $k$, the $k$-th tower is cut into $m_k$ subcolumns, $s_{k,r}$ spacer levels are placed above the $r$-th subcolumn, and the resulting stacks are placed from left to right. Since $\sum_{r=0}^{m_k-1}s_{k,r}=1$ for every $k$ and
\[
\sum_{k\ge0}(m_0m_1\cdots m_k)^{-1}<\infty,
\]
the total spacer measure is finite, so after normalization this yields a probability-preserving rank-one transformation. If $F_k$ denotes the base of the $k$-th tower and $h_k$ its height, then
\begin{equation}\label{eq:rank-one-height-recursion}
h_{k+1}=m_kh_k+\sum_{r=0}^{m_k-1}s_{k,r}=m_kh_k+1.
\end{equation}
\begin{lemma}\label{lem:rank-one-example-weak-mixing}
The transformation \(T\) constructed above is weak mixing.
\end{lemma}

\begin{proof}
Rank-one transformations are ergodic, so it suffices to rule out nontrivial eigenvalues. By the standard finite-measure rank-one eigenvalue criterion \cite[Section~4]{ChoksiNadkarni1995Eigenvalues}, an eigenvalue \(\lambda=e^{2\pi i\theta}\) must satisfy \(m_k\|h_k\theta\|_{\mathbb R/\mathbb Z}\to0\). Hence \(z_k:=\lambda^{h_k}\to1\) and \(z_k^{m_k}\to1\). Since \(h_{k+1}=m_kh_k+1\), we have \(\lambda^{h_{k+1}}=z_k^{m_k}\lambda\to\lambda\), while the same criterion gives \(\lambda^{h_{k+1}}\to1\). Thus \(\lambda=1\).
\end{proof}

For $D\in\mathscr B$ and $\tau>0$, let
\(
\mathcal E_{D,+}^T(\tau):=\{n\ge0:\ \mu(D\cap T^{-n}D)-\mu(D)^2>\tau\}.
\)
The following proposition is the main result of this subsection.

\begin{proposition}\label{prop:rank-one-polynomial-lower-bound}
For every \(\delta\in(0,1)\), there exist a set \(A\in\mathscr B\) and constants
\(C>0\), \(N_0\in\mathbb N\) such that every exceptional set \(J\) for \((A,A)\)
satisfies
\[
|J\cap[0,N]|\ge C N^{\delta}
\qquad(N\ge N_0).
\]
\end{proposition}

Fix $\delta\in(0,1)$, and let $M_j:=\lfloor m_j^{\delta}\rfloor$ and
\(
\rho_k:=\mu\Bigl(\bigcup_{r=0}^{h_k-1}T^rF_k\Bigr).\) Since the levels of the $k$-th tower exhaust $X$ modulo null sets, $\rho_k\to1$. Choose $K$ so large that
\begin{equation}\label{eq:rank-one-choice-K}
M_j\ge \frac12 m_j^{\delta}\quad (j\ge K),
\qquad
\sum_{j\ge K}\frac{2M_jh_j+2}{h_{j+1}}<\frac19,
\qquad
\rho_K>\frac67,
\qquad
h_K\ge5.
\end{equation}
Define \(\ell_K:=\lfloor h_K/2\rfloor\), $A:=\bigcup_{r=0}^{\ell_K-1}T^rF_K$, and
\[
\mathcal G:=\Bigl\{\sum_{j\ge K}d_jh_j:\ 0\le d_j<M_j\text{ for all }j,\ d_j=0\text{ for all but finitely many }j\Bigr\}.
\]
Then
\begin{equation}\label{eq:rank-one-alpha-gap}
\alpha:=\mu(A)=\rho_K\frac{\ell_K}{h_K}\in \left(\dfrac{1}{3}, \dfrac{1}{2}\right].
\end{equation}
\begin{lemma}\label{lem:rank-one-G-positive}
With the notation above, \(\mathcal G\subset \mathcal E_{A,+}^T(1/6).\)
\end{lemma}

\begin{proof}
For $j\ge K$ and $0\le d<M_j$, we claim that
\begin{equation}\label{eq:rank-one-shift-bound}
\mu(T^{dh_j}A\triangle A)\le \frac{2dh_j+2}{h_{j+1}}.
\end{equation}
For $q\ge j+1$, define recursively
\(
\Lambda_{j,j+1}:=\{0\}\) and \(\Lambda_{j,q+1}:=\bigsqcup_{a=0}^{m_q-1}(ah_q+\Lambda_{j,q})\). Fix $p>j+1$, let
\(
\mathcal T_p:=\bigcup_{r=0}^{h_p-1}T^rF_p,
\) and set
\[
\mathcal S_{j,p}:=\mathcal T_p\setminus \bigsqcup_{n\in\Lambda_{j,p}}\ \bigsqcup_{r=0}^{h_{j+1}-1}T^{n+r}F_p,
\qquad
\mathcal B(n,d):=T^{n+h_{j+1}-1}F_p\cup \bigsqcup_{c=0}^{d-1}\ \bigsqcup_{r=0}^{h_j-1}T^{n+ch_j+r}F_p,
\]
and
\begin{equation}
    \mathcal B_{j,p}(d):=\mathcal S_{j,p}\cup \bigcup_{n\in\Lambda_{j,p}}\mathcal B(n,d).
    \label{eqn0001}
\end{equation}
If $x\in \mathcal T_p\setminus \mathcal B_{j,p}(d)$, then $x\in T^{n+ch_j+r}F_p$ for some $n\in\Lambda_{j,p}$, $d\le c<m_j$, and $0\le r<h_j$, so
\[
T^{-dh_j}x\in T^{n+(c-d)h_j+r}F_p.
\]
Since $A$ is a union of levels of the $K$-th tower and $j\ge K$, membership in $A$ depends only on the relative level within a stage-$j$ tower copy; hence $1_A(x)=1_A(T^{-dh_j}x)$, and therefore
\begin{equation}
    (T^{dh_j}A\triangle A)\cap \mathcal T_p\subset \mathcal B_{j,p}(d).
    \label{eqn0002}
\end{equation}
Let $R_{j,p}$ be the number of levels of $\mathcal T_p$ lying in $\mathcal S_{j,p}$. Since the stage-$(j+1)$ copies indexed by $\Lambda_{j,p}$ are pairwise disjoint,
\begin{equation}
    \mu\Bigl(\bigcup_{n\in\Lambda_{j,p}}\mathcal B(n,d)\Bigr)
= |\Lambda_{j,p}|(dh_j+1)\frac{\rho_p}{h_p}
\le \frac{dh_j+1}{h_{j+1}}\,\rho_p.
\label{eqn0003}
\end{equation}
Note that $|\Lambda_{j,j+1}|=1$ and $R_{j,j+1}=0$, and for $q\ge j+1$,
\[
|\Lambda_{j,q+1}|=m_q|\Lambda_{j,q}|,
\qquad
R_{j,q+1}=m_qR_{j,q}+1.
\]
It follows that \(R_{j,p}\le |\Lambda_{j,p}|-1\). Since \(R_{j,p}\le |\Lambda_{j,p}|-1\) and \(|\Lambda_{j,p}|h_{j+1}\le h_p\), we have
\begin{equation}
    \mu(\mathcal S_{j,p})=R_{j,p}\frac{\rho_p}{h_p}
\le |\Lambda_{j,p}|\frac{\rho_p}{h_p}
\le \frac{\rho_p}{h_{j+1}}
\le \frac{dh_j+1}{h_{j+1}}\,\rho_p.
\label{eqn0004}
\end{equation}
By \eqref{eqn0001}, \eqref{eqn0002}, \eqref{eqn0003}, and \eqref{eqn0004}, we obtain
\[
\mu(T^{dh_j}A\triangle A) \leq \mu\bigl((T^{dh_j}A\triangle A)\cap \mathcal T_p\bigr)+1-\mu(\mathcal T_p)
\le \frac{2dh_j+2}{h_{j+1}}\rho_p+(1-\rho_p)\to \frac{2dh_j+2}{h_{j+1}}
\]
as $p\to\infty$, which proves \eqref{eq:rank-one-shift-bound}. For $n=\sum_{j\ge K}d_jh_j\in\mathcal G$, subadditivity, \eqref{eq:rank-one-shift-bound}, and \eqref{eq:rank-one-choice-K} give
\[
\mu(T^nA\triangle A)
\le \sum_{j:d_j>0}\mu(T^{d_jh_j}A\triangle A)
\le \sum_{j\ge K}\frac{2M_jh_j+2}{h_{j+1}}<\frac19.
\]
Using \eqref{eq:rank-one-alpha-gap},
\[
\mu(A\cap T^{-n}A)-\alpha^2
=\alpha-\alpha^2-\frac12\mu(T^nA\triangle A)
>\frac29-\frac1{18}=\frac16.
\]
Hence $n\in \mathcal E_{A,+}^T(1/6)$.
\end{proof}

\begin{lemma}\label{lem:rank-one-G-growth}
There exists $C>0$ such that
\[
|\mathcal G\cap[0,N]|\ge CN^{\delta}
\]  
for all sufficiently large $N$.
\end{lemma}

\begin{proof}
For $p>K$, let
\(
\mathcal G_p:=\Bigl\{\sum_{j=K}^{p-1}d_jh_j:\ 0\le d_j<M_j\Bigr\}.
\)
Since $M_j<m_j$, an induction using \eqref{eq:rank-one-height-recursion} gives
\(
\sum_{j=K}^{p-1}(M_j-1)h_j<h_p.
\)
Hence
\(
\mathcal G_p\subset[0,h_p)
\)
and
\(
\mathcal G_{p+1}=\bigsqcup_{t=0}^{M_p-1}(th_p+\mathcal G_p).
\)
We obtain
\begin{equation}\label{eq:rank-one-Gp-product-lower}
|\mathcal G_p|
=
\prod_{j=K}^{p-1}M_j
\ge
\Bigl(\prod_{j=K}^{p-1}m_j^{\delta}\Bigr)
\Bigl(\prod_{j=K}^{p-1}(1-m_j^{-\delta})\Bigr),\quad h_p^{\delta}
=
h_K^{\delta}
\Bigl(\prod_{j=K}^{p-1}m_j^{\delta}\Bigr)
\Bigl(\prod_{j=K}^{p-1}\Bigl(1+\frac1{m_jh_j}\Bigr)^{\delta}\Bigr).
\end{equation}
Since $\sum_{j\ge K}m_j^{-\delta}<\infty$ and $\sum_{j\ge K}(m_jh_j)^{-1}<\infty$, we see that
\(
\prod_{j=K}^{p-1}(1-m_j^{-\delta})\) and $\prod_{j=K}^{p-1}\Bigl(1+\frac1{m_jh_j}\Bigr)^{\delta}$ stay bounded away from $0$ and $\infty$, respectively. Therefore, by \eqref{eq:rank-one-Gp-product-lower}, there exists $c_1>0$ such that
\begin{equation}\label{eq:rank-one-Gm-growth}
|\mathcal G_p|\ge c_1h_p^{\delta}
\qquad (p>K).
\end{equation}
Now fix large $N$ and choose $p$ with $h_p\le N<h_{p+1}$. Let $n_p:=\lfloor N/h_p\rfloor$. Since $\mathcal G_p\subset[0,h_p)$, the sets $th_p+\mathcal G_p$ are pairwise disjoint for $0\le t<\min(n_p,M_p)$ and lie in $\mathcal G\cap[0,N]$, so
\begin{equation}\label{eqn:last}
|\mathcal G\cap[0,N]|
\ge \min(n_p,M_p)\,|\mathcal G_p|
\ge \frac12 n_p^\delta |\mathcal G_p|,
\end{equation}
where the second inequality follows because if \(n_p\le M_p\), then
\(n_p\ge n_p^\delta/2\), while if \(n_p>M_p\), then
\(n_p\le m_p\) since \(N<h_{p+1}=m_ph_p+1\), and hence \eqref{eq:rank-one-choice-K} yields
\(
M_p\ge  m_p^\delta/2\ge n_p^\delta/2.
\)
By \eqref{eq:rank-one-Gm-growth} and \eqref{eqn:last}, we obtain
\[
|\mathcal G\cap[0,N]|
\ge \frac12 n_p^\delta |\mathcal G_p|
\ge \frac12 c_1 (n_ph_p)^\delta
\ge \frac{1}{2^{1+\delta}}c_1N^\delta,
\]
where the last inequality follows from \(N<(n_p+1)h_p\le2n_ph_p\).
\end{proof}

\begin{proof}[Proof of Proposition~\ref{prop:rank-one-polynomial-lower-bound}]
By Lemmas~\ref{lem:rank-one-G-positive} and~\ref{lem:rank-one-G-growth}, there exist
\(c_1>0\) and \(N_1\in\mathbb N\) such that
\[
|\mathcal E_{A,+}^T(1/6)\cap[0,N]|\ge |\mathcal G\cap[0,N]|\ge c_1N^\delta
\]
for all \(N\ge N_1\). Let \(J\) be any exceptional set for \((A,A)\). Since
\(\mathcal E_{A,+}^T(1/6)\subseteq\mathcal E_A^T(1/6)\),
Lemma~\ref{lem:finite-modification-exceptional}(2) implies that
\[
\mathcal E_{A,+}^T(1/6)\setminus J
\]
is finite. Hence \(\mathcal G\setminus J\) is finite, and therefore
\[
|J\cap[0,N]|\ge C N^\delta
\]
for all sufficiently large \(N\), after decreasing the constant \(C>0\).
\end{proof}

\subsection{Additional Questions}
\label{subsec:open-questions}
We highlight some open problems related to our results.
\begin{enumerate}
 \item Can we generalize our results to all tight maps? For instance, is the condition $s_{m-1}=0$ necessary? We used this condition to ensure that the support of each $D_l$ is bounded. However, there are tight maps of interest, including Chacon's original construction of the Chacon map~\cite{chacon}, for which this condition fails. It is possible that a truncation argument could extend our results to this setting, but we leave this for future work.
 
 \item In Theorem~\ref{thm:lower-bound}, we showed the existence of a pair $(A,B)$ such that $J_{A,B}$ can be arbitrarily large. Does a contrasting phenomenon also occur for some pairs? More precisely, do there exist $A,B\in \mathscr B$ with $\mu(A),\mu(B)>0$ such that
 \[
 |J_{A,B}\cap[0,n]|\le C(\log n)^t
 \]
 for some $C,t>0$ and all sufficiently large $n$?
 
 \item For a given function $f(n)$, can we construct a weakly mixing map $T$ with
 \[
 C_1f(n)\le |J\cap[0,n]|\le C_2f(n)
 \]
 for some universal exceptional set $J$, or with
 \[
 C_1f(n)\le |J_{A,B}\cap[0,n]|\le C_2f(n)
 \]
 for some pair $(A,B)$?
 
 \item Beyond the class of tight maps, can one obtain sharper bounds for the size of exceptional sets in other weakly mixing systems? For interval exchange transformations, Corollary~\ref{cor:IET} gives an upper bound. Proposition~\ref{prop:rank-one-polynomial-lower-bound} shows that, in a weakly mixing one-spacer rank-one system, exceptional sets for a suitable pair can have polynomially large necessary growth. Can either phenomenon be sharpened? Can analogous lower-bound examples be established inside natural classes such as weakly mixing interval exchange transformations, or can stronger upper bounds be proved for broader classes of weakly mixing systems such as random substitution tilings, primitive substitution $\mathbb Z$-actions, and self-affine substitution tilings?
 
 \item Most of this paper is devoted to finding an exceptional set for a fixed weakly mixing transformation. The opposite question also seems natural: can one construct a weakly mixing system whose exceptional sets have prescribed optimal size? More specifically, for a given function $f(n)$, can one construct a weakly mixing map $T$ for which the optimal growth rate of exceptional sets is
 \[
 |J\cap[0,n]|\asymp f(n)?
 \]
 
 \item Is Proposition~\ref{prop:weak-conv-rate} optimal? More precisely, does there exist a weakly mixing system such that for every exceptional set $J_{A,B}$ and every sequence $c_n=o(nb_n)$, one has
 \[
 c_n\le |J_{A,B}\cap[0,n]|
 \]
 for all sufficiently large $n$?
\end{enumerate}

\appendix
\refstepcounter{section}
\section*{APPENDIX}

We define the total variation of a function and list some key properties.

\begin{definition}
 Let $f : \mathbb{R} \to \mathbb{R}$. We define the total variation $V(f)$ of $f$ as
 \[
 V(f) := \sup \left\{ \sum_{i=1}^{n-1} |f(x_i) - f(x_{i+1})| \right\},
 \]
 where the supremum is taken over all finite real numbers $x_1 < x_2 < \dots < x_n$. If $V(f) < \infty$, we say that $f$ has \emph{bounded variation}, and we denote the class of functions with bounded variation as $BV(\mathbb{R})$.
\end{definition}

We list some properties of the total variation.

\begin{lemma}
For any $f, g \in BV(\mathbb{R})$,
 \begin{gather}
 V(f \pm g) \le V(f) + V(g) \\
 V(\max(f, g)) \le V(f) + V(g),\qquad V(\min(f, g)) \le V(f) + V(g)
 \end{gather}
\end{lemma}

\begin{proof}
The first inequality follows from the triangle inequality on each partition.  For the second, use
\(\max(f,g)=(f+g+|f-g|)/2\) and \(V(|h|)\le V(h)\); the estimate for \(\min\) follows from \(\min(f,g)=-\max(-f,-g)\).
\end{proof}

\begin{lemma}
 Suppose $f \in BV(\mathbb{R}) \cap L^1(\mathbb{R})$ and $a>0$. Then,
 \[ 
 \left| \sum_{l \in \mathbb{Z}} f(x + al) - \frac{1}{a}\int_{-\infty}^{\infty} f(t)\,dt \right| \le V(f).
 \]
\end{lemma}

\begin{proof}
For \(I_l=[x+al,x+a(l+1)]\), put
\[
\epsilon_l=f(x+al)-\frac1a\int_{I_l}f(t)\,dt.
\]
Then
\[
|\epsilon_l|\le \frac1a\int_{I_l}|f(x+al)-f(t)|\,dt\le V(f;I_l).
\]
Summing over finite ranges and passing to the limit, using \(f\in L^1\), gives
\[
\left|\sum_l \epsilon_l\right|\le\sum_l V(f;I_l)\le V(f),
\]
which is the desired estimate.
\end{proof}

\begin{lemma}
Let \(\alpha\) be a probability distribution supported on a finite subset of
\(\mathbb R\), and let \(f\in BV(\mathbb R)\cap L^1(\mathbb R)\). If
\[
B(\alpha):=\inf\{R>0:\operatorname{supp}\alpha\subset[-R,R]\},
\]
then
\[
\|f-\alpha*f\|_1\le B(\alpha)V(f).
\]
\label{lem:conv-bound}
\end{lemma}

\begin{proof}
For \(\alpha=\delta_t\),
\[
\|f-\delta_t*f\|_1
=
\int_{\mathbb R}|f(x)-f(x-t)|\,dx
\le
|t|V(f).
\]
Indeed, this follows by decomposing \(\mathbb R\) into intervals of length
\(|t|\) and summing the variation along the corresponding arithmetic
progressions. For general finitely supported \(\alpha=\sum_i p_i\delta_{t_i}\), the triangle inequality gives
\[
\|f-\alpha*f\|_1
\le
\sum_i p_i\|f-\delta_{t_i}*f\|_1
\le
\sum_i p_i|t_i|V(f)
\le
B(\alpha)V(f).
\]
\end{proof}

\section*{Acknowledgments}
We would like to express our gratitude and thanks to Professor Seonhee Lim for her guidance throughout this project. We are also grateful to Professor Kyewon Koh Park for her insightful feedback. We would also like to thank Songun Lee for valuable discussions.

\bibliographystyle{plainurl}
\bibliography{ref}

\end{document}